\documentclass[11pt]{gtart}
\usepackage{amsmath, amssymb, graphicx, epsfig, graphs, verbatim, psfrag}

\newtheorem{thm}{Theorem}[section]

\newtheorem{cor}[thm]{Corollary}
\newtheorem{lem}[thm]{Lemma}
\newtheorem{prop}[thm]{Proposition}

\theoremstyle{definition}
\newtheorem{defn}[thm]{Definition}

\newtheorem{rem}[thm]{Remark}

\numberwithin{equation}{section}

\newcommand{\al}{\alpha}
\newcommand{\be}{\beta}
\newcommand{\ga}{\gamma}

\newcommand{\de}{\delta}

\newcommand{\om}{\omega}

\newcommand{\si}{\sigma}
\newcommand{\Si}{\Sigma}

\newcommand{\ze}{\zeta}
\newcommand{\bfz}{{\mathbb {Z}}}
\newcommand{\Ht}{{\mathcal {H}}}
\newcommand{\J}{{\mathcal {J}}}
\newcommand{\ot}{{\overline {\t}}}
\newcommand{\os}{{\overline {\s}}}

\newcommand{\x}{\times}
\newcommand{\s}{\mathbf s}

\renewcommand{\t}{\mathbf t}

\newcommand{\Z}{\mathbb Z}
\newcommand{\N}{\mathbb N}
\newcommand{\Q}{\mathbb Q}
\newcommand{\R}{\mathbb R}

\newcommand{\del}{\partial}
\newcommand{\hra}{\hookrightarrow}
\newcommand{\lra}{\longrightarrow}
\newcommand{\hf}{{{\widehat {HF}}}}

\newcommand{\Spin}{{\rm {Spin}}}

\begin{document}
\mathsurround=1pt 
\title{Tight contact structures on 
some small Seifert fibered 3--manifolds}

\author{Paolo Ghiggini} 
\address{CIRGET, Universit\'e du Qu\'ebec \`a Montr\'eal\\
Case Postale 8888, succursale Centre--Ville\\
Montr\'eal (Qu\'ebec) H3C 3P8, Canada} 
\email{ghiggini@math.uqam.ca}

\author{Paolo Lisca}\
\address{Dipartimento di Matematica ``L. Tonelli''\\
Largo Bruno Pontecorvo, 5\\
Universit{\`a} di Pisa \\I-56127 Pisa, Italy} 
\email{lisca@dm.unipi.it}

\author{Andr\'{a}s I. Stipsicz}
\address{R{\'e}nyi Institute of Mathematics\\
Hungarian Academy of Sciences\\
H-1053 Budapest\\ 
Re{\'a}ltanoda utca 13--15, Hungary}
\email{stipsicz@math-inst.hu}

\begin{abstract}
We classify tight contact structures on the small Seifert fibered
3--manifold $M(-1; r_1, r_2, r_3)$ with $r_i\in (0,1)\cap \Q$ and
$r_1, r_2\geq \frac{1}{2}$. The result is obtained by combining convex
surface theory with computations of contact Ozsv\'ath--Szab\'o
invariants.  We also show that some of the tight contact structures on the
manifolds considered are nonfillable, justifying the use of Heegaard
Floer theory.
\end{abstract}

\maketitle

\section{Introduction}
\begin{figure}[ht]
\begin{center}
\psfrag{e0}{\small $e_0$}
\psfrag{r1}{$-\frac{1}{r_1}$}
\psfrag{r2}{$-\frac{1}{r_2}$}
\psfrag{r3}{$-\frac{1}{r_3}$}
\includegraphics[height=2cm]{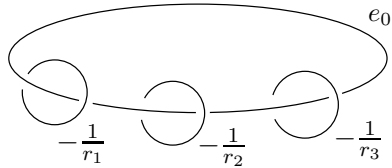}
\end{center}
\caption{\quad Surgery diagram for the Seifert 
fibered 3--manifold $M(e_0; r_1, r_2, r_3)$}
\label{f:seifert}
\end{figure}

In this paper we call a Seifert fibered 3--manifold $M$~\emph{small}
if $M$ is closed, the base surface is $S^2$ and $M$ has exactly three
singular fibers. Using normalized Seifert invariants, a small Seifert
fibered 3--manifold $M$ can be described by the surgery diagram of
Figure~\ref{f:seifert}, where $e_0\in \bfz $ and $r_i \in (0,1)\cap
\Q$, with $r_1\geq r_2\geq r_3$. Conversely, a 3--manifold given by
Figure~\ref{f:seifert} carries a natural structure of a small  Seifert
fibered 3--manifold. We shall denote such a 3--manifold by $M(e_0;
r_1, r_2, r_3)$.
The classification of tight contact structures on $M(e_0;r_1,r_2,r_3)$
has been given in~\cite{Wu2} when $e_0\neq 0, -1, -2$, and then
extended in~\cite{gls} to the case $e_0\geq 0$. In the present article we
consider a further special case: 
\[
e_0=-1,\quad r_1\geq r_2\geq\frac{1}{2}.  
\]
Let us assume that $-\frac{1}{r_i}$ has continued fraction
expansion $[a_0^i, \ldots , a_{k_i}^i]$, that is,
\[ 
-\frac{1}{r_i}=-a_0^i-\cfrac1{-a_1^i-\cfrac1{\ddots - \cfrac1{-a_{k_i}^i}}},
\quad
a^i_j\in\N,\ a^i_j\geq 2.
\]

Define 
\[
h(r_1,r_2,r_3):= (a_1^3-1)\prod_{i=1}^3 \prod_{j\geq 2}(a_j^i-1),
\]
\[
\varphi(r_1,r_2,r_3):=
\begin{cases}
\left(2(a_1^1-1)(a_1^2-1) + (a_0^3-1)(a_1^1+a_1^2-2)\right)h(r_1,r_2,r_3)
\quad\text{if}\quad r_2>\frac 12,\\
\left(2(a_1^1 -1) + (a_0^3 -1)\right)(a_1^3-1)\prod_{i\neq2}
\prod_{j\geq 2}(a^i_j-1)
\quad\text{if}\quad r_1>r_2=\frac 12,\\
2\prod_{j\geq 1}(a_j^3-1)\quad\text{if}\quad r_1=r_2=\frac12,
\end{cases}
\]
and 
\[
\psi(r_1,r_2,r_3):=
\begin{cases}
(a_1^1-1)(a_1^2-1)a_1^3\prod_{i=1}^3\prod_{j\geq 2}(a_j^i-1)
\quad\text{if}\quad r_3\neq\frac1{a_0^3},\\
\prod_{i=1}^2\prod_{j\geq 1}(a_j^i-1)
\quad\text{if}\quad r_3 = \frac1{a_0^3},
\end{cases}
\]
where one should conventionally set $a^i_j=2$ when $j>k_i$.
Our main result is the following.

\begin{thm}\label{t:main}
Let $M=M(-1; r_1, r_2, r_3)$ be a small Seifert fibered 3--manifold with
$r_1\geq r_2\geq\frac12$. Then, $M$ supports exactly  
\[ 
\varphi(r_1,r_2,r_3) + \psi(r_1,r_2,r_3)
\]
tight contact strucures up to isotopy.
\end{thm}

In order to put this result in perspective, recall that the classification of
tight contact structures on small Seifert fibered 3--manifolds with $e_0\neq
-1, -2$ (as it is given in \cite{gls, Wu2}) required two steps. First, using
convex surface theory, an upper bound on the number of tight contact
structures was achieved.  In the second step, using Legendrianizations of
appropriate surgery diagrams for the 3--manifold at hand, a collection of
Stein fillable contact structures was constructed. The isotopy classes of
these Stein fillable (hence tight) structures were distingushed by the first
Chern classes of their Stein fillings, resting on a result of \cite{LM}.  For
the class of 3--manifolds considered in this paper, this strategy needs to be
modified, since --- as it will be shown in Section~\ref{s:special},
cf.~Theorem~\ref{t:nonfill} --- some 3--manifolds addressed by
Theorem~\ref{t:main} admit tight, nonfillable contact structures. Tightness of
such nonfillable structures can be either shown by state traversal methods
(which becomes extremally complicated when the underlying 3--manifold is
atoroidal, as is the case of the 3--manifolds considered in
Theorem~\ref{t:main}), or by computing contact Ozsv\'ath--Szab\'o
invariants. By finding appropriate contact $(\pm 1)$--surgery diagrams for all
the potential tight contact structures, here we describe a simple way to
determine their contact Ozsv\'ath--Szab\'o invariants.

Our results imply that the mod 2 reduced contact Ozsv\'ath--Szab\'o
invariant is a complete invariant for tight contact structures on the
3--manifolds $M(-1;r_1,r_2,r_3)$, $r_1,r_2\geq\frac12$. Moreover, two
tight structures on these 3--manifolds are isotopic if and only if
they induce the same spin$^c$ structure.

Notice that the class of small Seifert fibered 3--manifolds considered
in Theorem~\ref{t:main} is large enough to contain each small Seifert
fibered 3--manifold with finite fundamental group (with one
orientation) --- with the unique exception of the Poincar\'e homology
sphere $-\Sigma (2,3,5)=M(-1; \frac12, \frac13, \frac15)$.

In Section~\ref{s:convex} we derive upper bounds for the number of
isotopy classes of tight contact structures on the Seifert fibered
3--manifolds listed in Theorem~\ref{t:main}. The heart of these
arguments is to establish isotopies between tight contact structures
presented in different ways. Section~\ref{s:OSz} is devoted to a
recollection of relevant results of Heegaard Floer theory related to
the contact Ozsv\'ath--Szab\'o invariants. Contact surgeries are also
briefly discussed.  In Section~\ref{s:special} the special case of
$M(-1; \frac{1}{2}, \frac{1}{2}, \frac{1}{p})$ is examined, and by
using contact Ozsv\'ath--Szab\'o invariants we achieve a complete
classification for these Seifert fibered 3--manifolds. Stein
fillability and nonfillability is also discussed. Finally in
Section~\ref{s:general} we complete the proof of Theorem~\ref{t:main}
by giving lower bounds for the number of tight contact structures on
$M(-1; r_1, r_2, r_3)$ satisfying $r_1, r_2\geq \frac{1}{2}$.

There is a final remark in place. Notice that in the definition of
continued fractions we introduced negative signs, so that the
coefficients $a^i_j$ became positive. Although this convention might
differ from many results established in the literature, for our
purposes it seemed to be more convenient to work with positive rather
than negative numbers.

{\bf Acknowledgements.} Part of this work was done while the first
author was visiting Princeton University supported by the NSF Focused
Research Grant FRG-024466.  The third author wants to thank the
Institute for Advanced Study, Princeton for their hospitality while
part of this collaboration was carried out. He was partially supported
by NSF Focused Research Grant FRG-024466 and by OTKA T049449.

\section{Upper bounds from convex surface theory}
\label{s:convex}

In order to get an upper bound for the number of tight contact
structures, we will follow the methods developed in~\cite{EH, Gi0} and
implemented~\cite{EH, GhS}.  Suppose that $(M,\xi)$ is a contact,
Seifert fibered $3$--manifold. Then, a Legendrian knot in $M$ smoothly
isotopic to a regular fiber admits two framings: one coming from the
fibration and another one coming from the contact structure $\xi$. The
difference between the contact framing and the fibration framing is
the \emph{twisting number} of the Legendrian curve. We say that $\xi$
has \emph{maximal twisting equal to zero} if there is a Legendrian
knot $L$ isotopic to a regular fiber such that $L$ has twisting number
zero. Applying~\cite[Theorem~1.5(1)]{Wu1} we have:

\begin{prop}\label{p:wu}
Every tight contact structure on $M(-1; r_1 , r_2, r_3)$ with $r_1\geq
r_2\geq \frac{1}{2}$ has maximal twisting equal to zero.
\end{prop}

\begin{proof}
Using unnormalised Seifert invariants, we write $M(-1; r_1, r_2,
r_3)$ as $M(0; r_1- 1, r_2, r_3)$.  In the notation
of~\cite[Theorem~1.5]{Wu1} this means that $\frac{q_3}{p_3}=r_1-1$ and
$\frac{q_2}{p_2}=r_2$. Therefore
$\frac{q_2}{p_2}+\frac{q_3}{p_3}=r_1+r_2-1$, which is nonnegative by
our assumptions. Consequently~\cite[Theorem~1.5(1)]{Wu1} applies,
implying the statement.
\end{proof}

Let $F_i$ ($i=1,2,3$) be the three singular fibers of the Seifert
fibration on $M$. First, in view of Proposition~\ref{p:wu} we can 
isotope the Seifert fibration until there is a Legendrian regular fiber 
$L$ with twisting number zero with respect to the fibration.
Then, we can isotope each $F_i$ further so that it becomes
Legendrian. 

Let $V_i$ be a standard convex neighbourhood of $F_i$,
$i=1,2,3$. Then, $M \setminus (V_1 \cup V_2 \cup V_3)$ can be
identified with $\Sigma \times S^1$ where $\Sigma$ is a
pair of pants. This diffeomorphism determines identifications of $-
\partial (M \setminus V_i)$ with $\R^2 / \Z^2$ so that $\binom{1}{0}$
is the direction of the section $\Sigma \times \{ 1 \}$ and
$\binom{0}{1}$ is the direction of the regular fibers. In order to fix
one among the infinitely many product structures on $M \setminus (V_1
\cup V_2 \cup V_3)$ we also require the meridian of each $V_i$ to have
slope $- \frac{\beta_i}{\alpha_i}$ in $- \partial (M \setminus V_i)$,
with
\[
\frac{\beta_1}{\alpha_1} =r_1-1,\quad
\frac{\beta_2}{\alpha_2} = r_2,\quad
\frac{\beta_3}{\alpha_3} = r_3.
\]
We also choose an identification between $\partial V_i$ and $\R^2 /
\Z^2$ so that $\binom{1}{0}$ is the direction of the meridian of $V_i$
and $\binom{1}{0}$ is the direction of a longitude. Notice that
$\partial V_i$ and $- \partial (M \setminus V_i)$ coincide as sets,
but are identified with $\R^2 / \Z^2$ in different ways. We can choose
the longitude on $V_i$ so that these two identifications are
related by gluing matrices $A_i: \partial V_i \to - \partial (M
\setminus V_i)$ given by
\[
A_i= \left ( \begin{matrix}
\alpha_i & \alpha_i' \\ 
- \beta_i & -\beta_i'  
\end{matrix} \right )
\]
with $\beta_i \alpha_i'-\alpha_i \beta_i'=1$ and $0< \alpha_i' < \alpha_i$.

Since $V_i$ is a standard convex neighbourhood of a Legendrian curve,
$- \partial (M \setminus V_i)$ is a standard torus with exactly $2$
dividing curves. By flexibility of the Legendrian
ruling~\cite[Corollary~3.6]{H1} we can modify the characteristic
foliation of $- \partial (M \setminus V_i)$ so that its Legendrian
ruling has infinite slope.  Consider convex vertical annuli $A_i$ with
Legendrian boundary between a vertical Legendrian ruling curve of $-
\partial (M \setminus V_i)$ and the Legendrian regular fiber $L$ with
twisting number zero. By the Imbalance Principle
\cite[Proposition~3.17]{H1} the annuli $A_i$ give bypasses on $-
\partial (M \setminus V_i)$. The Bypass Attachment Lemma
\cite[Lemmas~3.12 and~3.15]{H1} implies that the attachments of these
bypasses reduce the denominator of the slope $- \partial (M \setminus
V_i)$ in absolute value, therefore after a finite number of bypass
attachments we obtain tubular neighbourhoods $U_i$ of $F_i$ containing
$V_i$ such that $- \partial (M \setminus U_i)$ has infinite slope for
$i=1,2,3$.

\begin{lem}\label{l:slopeUi}
The convex torus $\partial U_i$ has slope 
\[
- \frac{\alpha_i}{\alpha_i'}= [a^i_{k_i}, \ldots , a^i_0],
\quad i=1,2,3.
\]
\end{lem}

\begin{proof}
For $i=2,3$ the statement follows from \cite[Lemma~A4]{OW} because
$r_i=\frac{\beta_i}{\alpha_i}$. When $i=1$ we can still
apply~\cite[Lemma~A4]{OW} because
$r_i=\frac{\beta_1-\alpha_1}{\alpha_1}$ and we have 
\[
\alpha_1'(\beta_1 - \alpha_1)-\alpha_1 (\beta_1'-\alpha_1')=1.
\]
\end{proof}

Since $- \frac{\al_i}{\al_i'} < -1$, by~\cite[Proposition~4.16]{H1}
 applied to $U_i \setminus V_i$ there exist tubular neighbourhoods
 $V_i'\subset U_i$ of the singular fibers such that $\partial V_i'$ is
 a standard convex torus with slope
 $-1$. Following~\cite[Section~4.4]{H1} we decompose $U_i \setminus V'_i$ into 
consecutive layers $N_j^i$ diffeomorphic to toric annuli
 with convex boundary and boundary tori with slopes
\[
[a^i_{k_i}, \ldots, a^i_1-1]\quad\text{and}\quad 
[a^i_{k_i}, \ldots , a^i_0]\quad\text{if}\quad j=0, 
\]
\[
[a^i_{k_i}, \ldots , a^i_{j+1}-1]\quad\text{and}\quad
[a^i_{k_i}, \ldots , a^i_j-1]\quad\text{if}\quad 0 < j < k_i, 
\]
\[
-1 \quad \text{and} \quad [a^i_{k_i}-1] \quad \text{if} \quad j=k_i, 
\]
where $[a^i_{k_i}, \ldots, a^i_1-1]$ should be interpreted as $-1$ when 
$r_i=\frac1{a_0^i}$. If $a^i_j = 2$ then by definition $N_j^i$ is an 
invariant neighbourhood of a convex torus, and if $a^i_j > 2$ 
by~\cite[Proposition~4.14]{H1} $N^i_j$ is a continued fraction block.

We define $q^i_j$ as the number of positive basic slices in
$N^i_j$. We have 
\[
0 \leq q^i_0 \leq a^i_0-1\quad\text{and}\quad 
0 \leq q^i_j \leq a^i_j -2  \qquad (j>0) .  
\]
For $j > k_i$ there are no more layers $N_j^i$ and we define $q^i_j=
\infty$.

Our assumption $r_1 \geq r_2 \geq \frac 12$ implies that $a_0^1 =
a_0^2 = 2$, therefore $q^1_0, q^2_0 \in \{ 0,1\}$ and $N_0^1$,
$N_0^2$ are basic slices.

By the classification of tight contact structures on solid
tori~\cite[Theorem~2.3]{H1} there are 
\[
a_0^i(a_1^i -1) \ldots (a^i_{k_i} -1) 
\]
distinct tight contact structures on each $U_i$. We will see that not
every combination of tight contact structures on $U_1$, $U_2$, and
$U_3$ is the restriction of a tight contact structure on $M$. But
first we need to analyse the tight contact structures on the
complement of $U_1$, $U_2$, and $U_3$. There are infinitely many
distinct tight contact structures on $M \setminus (U_1 \cup U_2 \cup
U_3)$ which can be the restriction of a tight contact structure on
$M$. However, Lemma~\ref{l:struttura} below implies that the isotopy
class of a tight contact structure $\xi$ on $M$ does not depend on its
restriction to $M \setminus (U_1 \cup U_2 \cup U_3)$.
 
\begin{defn}
Let $\Sigma$ be a pair of pants. A tight contact structure $\xi$ on
$\Sigma \times S^1$ is {\em appropriate} if there is no contact
embedding 
\[
(T^2 \times I,\xi_{\pi})\hra(\Sigma\x S^1,\xi)
\]
with $T^2 \times \{0\}$ isotopic to a boundary component, where
$\xi_{\pi}$ is a tight contact structure with convex boundary and
twisting $\pi$ (see~\cite[{\S}~2.2.1]{H1} for the definition of
twisting).
\end{defn}

\begin{lem}\label{l:struttura}
Let $\Sigma$ be a pair of pants and let $\xi$ be an appropriate
contact structure on $\Sigma \times S^1$ with convex boundary $-
\partial (\Sigma \times S^1)= T_1 \cup T_2 \cup T_3$, boundary slopes
\[
s(T_1)=-n,\quad s(T_2)= -1,\quad 
s(T_3)=\infty,\quad n \in \N \cup \{ 0 \},
\]
and $\#\Gamma_{T_i}=2$ for $i=1,2,3$. Then, there is a pair of pants
$\Sigma'$ contained in $\Sigma$ and a factorization 
\[
\Sigma \times S^1= (\Sigma' \times S^1) \cup B_1 \cup B_2 
\]
such that
\begin{enumerate}
\item 
the restriction of $\xi$ to $\Sigma'\times S^1$ is appropriate and has
convex boundary with infinite boundary slopes;
\item 
the restrictions of $\xi$ to $B_1$ and $B_2$ are basic slices;
\item 
the isotopy class of $\xi$ is determined by the signs of the restrictions 
of $\xi$ to $B_1$ and $B_2$.
\end{enumerate}
\end{lem}
\begin{figure}[htb]
\begin{center}
\includegraphics[width=12cm]{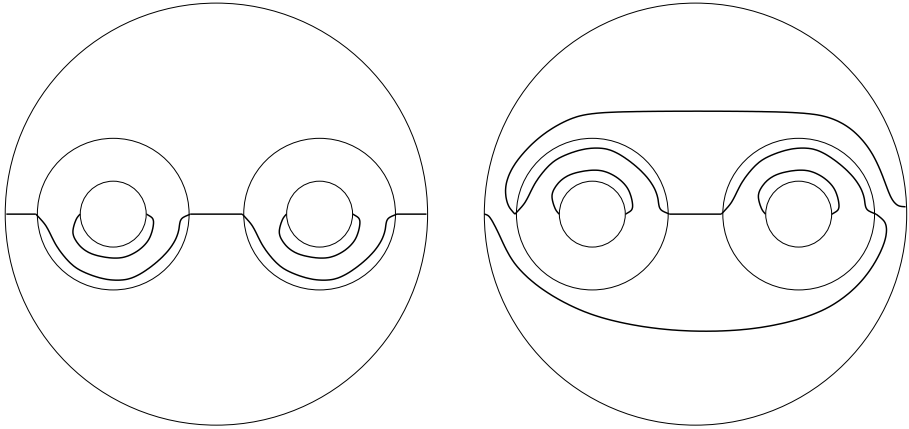}
\end{center}
\caption{Nonisotopic dividing sets on a convex horizontal section of
$\Sigma' \times S^1$ become isotopic when extended to a section
of $\Sigma'' \times S^1$}
\label{f:curves}
\end{figure}
\begin{proof}
The proof of this lemma is similar to the proof
of~\cite[Lemma~4.1]{Wu2}, where different boundary slopes are
considered.  Here we give a sketch of the argument.

The existence of the factorization is stated
in~\cite[Lemma~5.1(a)]{H2}. The restriction of $\xi$ to $\Sigma'
\times S^1$ is appropriate because the restriction of an appropriate
contact structure is always appropriate. This proves (1). $B_1$ and
$B_2$ are basic slices because they are appropriate (which, for
thickened tori, is the same as being minimally twisting) and their
boundary slopes are consecutive in the Farey Tessellation. This proves
(2). Let $\Sigma_0'$ be a convex horizontal section of $\Sigma' \times
S^1$ with Legendrian boundary consisting of horizontal Legendrian
ruling curves of $T_1$, $T_2$, and $T_3$ so that the dividing set
$\Gamma_{\Sigma_0'}$ has the minimum possible number of dividing
curves in the isotopy class of $\Sigma'_0$. The dividing set of
$\Sigma_0'$ cannot contain a boundary parallel dividing curve because
$\xi$ is appropriate, therefore it consists of three arcs, each
joining distinct boundary components: see~\cite[Lemma~10]{EH}.
By~\cite[Lemma~5.2]{H2} there is a unique extension of $(\Sigma \times
S^1, \xi)$ to $(\Sigma'' \times S^1, \xi'')$ obtained by adding basic
slices $B_1''$ and $B_2''$ so that the resulting contact structure
$\xi''$ is tight and has infinite boundary slopes. Let $\Sigma_0''$ be
a convex horizontal section of $\Sigma'' \times S^1$ with Legendrian
boundary extending $\Sigma_0'$ so that the dividing set
$\Gamma_{\Sigma_0''}$ has the minimum number of dividing curves in the
isotopy class of $\Sigma_0''$. Then the dividing set
$\Gamma_{\Sigma_0''}$ consists of three boundary parallel arcs, as in
Figure~\ref{f:curves}. Therefore, the isotopy class of $\xi''$ is
determined by the signs of the boundary parallel regions. These signs
depend on the signs of the basic slices $B_1''$ and $B_2''$, which
must be equal to the signs of the basic slices $B_1$ and $B_2$ by the
Gluing Theorem~\cite[Theorem~4.25]{H1}.
\end{proof}

\begin{lem} \label{grafico}
Let $(T^2 \times [0,1], \eta)$ be a minimally twisting tight contact
structure with boundary slopes $s_0$ and $s_1$. Let $T^2 \times [\frac
12, 1]$ be the last basic slice in the basic slice decomposition of
$(T^2 \times [0,1], \eta)$. Then, the slope $s_{\frac12}$ of $T^2
\times \{ \frac 12 \}$ is the vertex of the Farey Tessellation in the
counterclockwise arc starting at $s_{\frac12}$ and ending at $s_1$
which is the closest to $s_0$ among those connected to $s_1$ by an
edge.
\end{lem}

\begin{proof}
In \cite[Section~4.4.3]{H1} $T^2 \times \{ \frac 12 \}$ is obtained from
$T^2 \times \{ 1 \}$ by attaching a bypass along a Legendrian ruling curve
with slope $s_1$. The bypass attachment lemma~\cite[Lemma~3.15]{H1}
concludes the proof.
\end{proof}

Let $Z_i$ be the complement of the outermost basic slice in $U_i$. In
particular, $U_1\setminus Z_1$ and $U_2 \setminus Z_2$ coincide with
$N_0^1$ and $N_0^2$ respectively. 

\begin{lem}\label{l:Zi}
$- \partial (M \setminus Z_1)$ has boundary slope $0$, while $-
\partial (M \setminus Z_2)$ and $- \partial (M \setminus Z_3)$ have
both boundary slopes equal to $-1$.
\end{lem}

\begin{proof}
Since
\[
A_1
\begin{pmatrix}
\phantom{-}1\\
-1
\end{pmatrix}=
\begin{pmatrix}
\al_1-\al'_1\\
-\be_1+\be'_1
\end{pmatrix},
\]
using the relation $\al'_1\be_1-\al_1\be'_1=1$ and the fact that 
\[
\frac{\be_1}{\al_1} = r_1 -1 < 0\quad\text{implies}\quad
\be_1\leq -1,
\]
we see that the slope of $-\partial (M \setminus V'_i)$ is
\[
\frac{-\be_1+\be'_1}{\al_1-\al'_1}=
\frac{-\al_1\be_1+\al_1\be'_1}{\al_1(\al_1-\al'_1)}=
\frac{-\be_1(\al_1-\al'_1)-1}{\al_1(\al_1-\al'_1)}=
-\frac{\be_1}{\al_1}-\frac1{\al_1(\al_1-\al'_1)} \geq 
-\frac{\be_1+1}{\al_1}\geq 0.
\] 
On the other hand, 
\[
\frac{-\be_1+\be'_1}{\al_1-\al'_1}= 
-\frac{\be_1}{\al_1}-\frac1{\al_1(\al_1-\al'_1)}
\leq -\frac{\be_1}{\al_1} = 1-r_1 < 1.
\]
Applying Lemma~\ref{grafico} to $U_1\setminus V_1'$ with its boundary
slopes computed with respect to the basis of $- \partial (M \setminus
U_1)$, we find that $- \partial (M \setminus Z_1)$ has boundary slope
$0$. In the same way we see that $- \partial (M \setminus Z_2)$ and $-
\partial (M \setminus Z_3)$ have both boundary slope $-1$.
\end{proof}

It is easy to check that if $\xi$ is a tight contact structure on $M$,
its restriction to $M \setminus (Z_1 \cup Z_2 \cup U_3)$ is
appropriate, therefore Lemma~\ref{l:struttura} implies that the
isotopy class of a tight contact structure $\xi$ on $M$ 
depends only on the
restriction of $\xi$ to the solid tori $U_i$. This fact implies that
the number of possible tight contact structures on $M$ is at most 
\[
\prod_{i=1}^3 a_0^i \prod_{j\geq 1}(a_j^i-1).
\]
We will sharpen this upper bound by excluding some overtwisted contact
structures and eliminating the overcounting of different presentations
of the same tight contact structure. Lemmas \ref{l:backward}
and  \ref{l:scambio} are essentially  more precise reformulations of
a particular case of \cite[Theorem 4.13]{GhS}.  

\begin{lem}\label{l:backward}
Let $\Sigma$ be a pair of pants and let $\xi$ be a
contact structure on $\Sigma \times S^1$ with convex boundary $-
\partial (\Sigma \times S^1)= T_1 \cup T_2 \cup T_3$, boundary slopes
\[
s(T_1)=-n,\quad s(T_2)= -1,\quad 
s(T_3)=\infty,\quad n \in \N \cup \{ 0 \},
\]
and $\#\Gamma_{T_i}=2$ for $i=1,2,3$. Suppose that there exists a
collar neighbourhood $B_3\subset \Sigma \times S^1$ of $T_3$ such that
\begin{enumerate}
\item 
$B_3$ is a basic slice with boundary slopes $\infty$ and $n$; 
\item 
the restriction of $\xi$ to $(\Sigma \times S^1) \setminus B_3 =
\Sigma'' \times S^1$ coincides, up to isotopy, with the unique tight
contact structure on $\Sigma'' \times S^1$ without vertical Legendrian
curves with twisting number $0$, and with boundary slopes $-n$, $-1$,
and $n$.
\end{enumerate}
Then $(\Sigma \times S^1, \xi)$ is appropriate, and in the
decomposition $\Sigma \times S^1= B_1 \cup B_2 \cup (\Sigma' \times
S^1)$ of Lemma \ref{l:struttura}, where $B_1$ and $B_2$ are basic
slices with boundary slopes $-n$, $\infty$ and $-1$, $\infty$
respectively, the signs of $B_1$ and $B_2$ are both opposite to the
sign of $B_3$.
\end{lem} 

\begin{figure}[ht]
\begin{center}
\psfrag{1}{$B_1$}
\psfrag{2}{$B_2$}
\psfrag{3}{$B_3$}
\includegraphics[width=12cm]{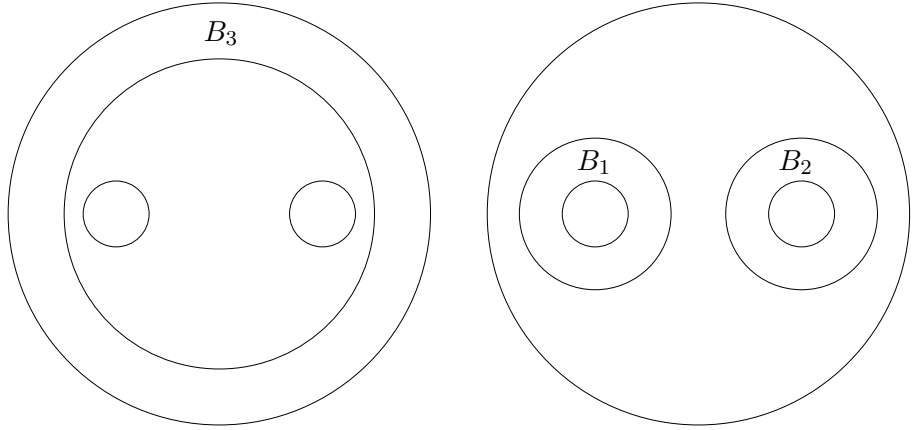}
\end{center}
\caption{The two decompositions of $(\Sigma \times S^1, \xi)$.}
\label{f:curves2}
\end{figure}

Here we warn the reader that our convention for the computation of the
boundary slopes differs from that of~\cite[Lemma~5.1]{H2} in the sense
that our slopes are computed with respect to $- \partial (\Sigma
\times S^1)$, as opposed to $\partial (\Sigma \times S^1)$.

\begin{proof}
Define $- \partial B_i= T_i - T_i'$ for $i=1,2,3$.  $(\Sigma'' \times
S^1, \xi|_{\Sigma'' \times S^1})$ is contactomorphic to the complement
of a vertical Legendrian curve with twisting number $-1$ in a
nonrotative thickened torus with boundary slopes $-n$ because there is
a unique tight contact structure up to isotopy on $\Sigma'' \times
S^1$ with these boundary slopes and without vertical Legendrian curves
with twisting number $0$, see \cite[Lemma~5.1(4b)]{H2}. A collar of
$T_3'$ in $B_3$ is isomorphic to a nonrotative thickened torus,
therefore we can identify $(\Sigma \times S^1, \xi)$ with the complement
of a vertical Legendrian curve with twisting number $-1$ in
$B_3$. Since $B_3$ is tight and minimally twisting we can conclude
that $(\Sigma \times S^1, \xi)$ is appropriate.

To prove that the signs of the restrictions of $\xi $ to $B_1$ and
$B_2$ are the opposite of the sign of its restriction to $B_3$, we
argue by evaluating the relative Euler class $e(\xi)$ of $\xi$ on
vertical annuli $A_1 \subset B_1$ and $A_2 \subset B_2$ with
Legendrian boundary
\[
\partial A_i= (A_i \cap T_i') - (A_i \cap T_i)
\]
such that $A_i \cap T_i$ is a Legendrian ruling curve of $T_i$, and
$A_i \cap T_i'$ is a Legendrian divide of $T_i'$.
 
Let $A_i'$, for $i=1,2$, be a vertical annulus in $\Sigma' \times S^1$
such that $\del A_i'= (A_i' \cap T_3)- (A_i' \cap T_i')$ consists of
Legendrian divides of $T_i'$ and $T_3$, and $A_i' \cap T_i' = A_i \cap
T_i'$. We consider also the vertical annuli $A_i''=A_i \cup A_i'$
between $T_3$ and $T_i$.  From the Thurston--Bennequin inequality it
follows that
\[
\langle e(\xi),[A_i']\rangle = 0,
\]
therefore
\[
\langle e(\xi), [A_i] \rangle = \langle e(\xi), [A_i''] \rangle.
\]
Take a vertical annulus $A_3 \subset B_3$ with Legendrian boundary
$\partial A_3= (A_3 \cap T_3')-(A_3 \cap T_3)$, and a vertical annulus
$A_3' \subset \Sigma'' \times S^1$ with Legendrian boundary $\partial
A_3'= (A_3 \cap T_i)-(A_3 \cap T_3')$ for either $i=1,2$. Assume
moreover that $A_3 \cap T_3'=A_3' \cap T_3'$. The dividing set of
$A_3'$ can contain no boundary parallel dividing arc because $\Sigma''
\times S^1$ contains no vertical Legendrian curve with twisting number
$0$, therefore \cite[Proposition~4.5]{H1} implies that $\langle
e(\xi), [A_3'] \rangle=0$.  If we call $A_3''=A_3 \cup A_3'$, then we
have
\[
\langle e(\xi), [A_3] \rangle = \langle e(\xi), [A_3''] \rangle = -
\langle e(\xi), [A_i''] \rangle.
\]
The change of sign in the evaluation of the relative Euler class is
due to the different orientations of $A_3''$ on one hand, and $A_i''$
with $i=1,2$ on the other hand.  This implies that
\[
\langle e(\xi), [A_3] \rangle = - \langle e(\xi), [A_i] \rangle
\]
for $i=1,2$.
\end{proof}

Combining Lemma \ref{l:struttura} and Lemma \ref{l:backward} we obtain
the following corollary, which contains the basic move which allows us
to go from a sign configuration of the basic slice decompositions of
$U_1$, $U_2$, and $U_3$ to a different one without affecting the
isotopy class of $\xi$.

\begin{lem}\label{l:scambio}
Let $\Sigma$ be a pair of pants and let $\xi$ be an appropriate
contact structure on $\Sigma \times S^1$ with convex boundary 
$- \partial (\Sigma \times S^1)= T_1 \cup T_2 \cup T_3$, boundary slopes
\[
s(T_1)=-n,\quad s(T_2)= -1,\quad s(T_3)=\infty,\quad n \in \N \cup \{
0 \},
\]
and $\#\Gamma_{T_i}=2$ for $i=1,2,3$. If the basic slices $B_1$ and
$B_2$ of Lemma~\ref{l:struttura} have the same sign, then there exists
a collar neighbourhood $B_3$ of $T_3$ such that $B_3$ is a basic slice
with boundary slopes $\infty$ and $n$ having sign opposite to the sign
of $B_1$ and $B_2$ and the restriction of $\xi$ to $\Sigma'' \times
S^1 = (\Sigma \times S^1) \setminus B_3$ coincides, up to isotopy,
with the unique tight contact structure on $\Sigma'' \times S^1$
without vertical Legendrian curves with twisting number $0$, and with
boundary slopes $-n$, $-1$, and $n$. 
\end{lem}

\begin{proof}
Let $\Si''$ be a pair of pants, and let $\eta'$ denote the the unique
tight contact structure on $\Si''\x S^1$ without vertical Legendrian
curves with twisting number $0$, and with boundary slopes $-n$, $-1$,
and $n$. Consider the contact structure $\eta$ obtained by gluing to
the boundary component of $(\Si''\x S^1,\eta')$ with slope $n$ a basic
slice $B_3$ with boundary slopes $\infty$ and $n$, and with sign
opposite to the sign of $(B_1,\xi|_{B_1})$ and
$(B_2,\xi|_{B_2})$. Then, by Lemma~\ref{l:backward}, in the
decomposition of Lemma~\ref{l:struttura} for $\eta$ the signs of the
basic slices $B_1$ and $B_2$ are both opposite to the sign of $B_3$,
and therefore are equal to the sign of $(B_1,\xi|_{B_1})$ and
$(B_2,\xi|_{B_2})$. By Lemma~\ref{l:struttura} the contact structures
$\eta$ and $\xi$ are isotopic, so the statement follows.
\end{proof}

We will see that the nontrivial behaviour of the tight contact
structures on $M$ comes from the first two outermost layers $N_0^i$
and $N_1^i$ in the decomposition of $U_i$ ($i=1,2,3$). This fact
justifies the introduction of the following notation: to a potentially
tight contact structure on $M$ we associate the matrix
\[
\left (
\begin{matrix} q_0^1 & q_0^2 & q_0^3 \\ q_1^1 & q_1^2 & q_1^3
\end{matrix} \right )
\]
whose entries $q^i_j$ are the number of positive basic slices in the
basic slice decompositions of $N_j^i$ for $i=1,2,3$ and
$j=0,1$. (Recall that we defined $q_j^i=\infty$ if $j>k_i$.) We will
call this matrix {\em the matrix of signs} of the contact structure.

We will study two separate cases. Suppose first that $q^1_0 \neq
q^2_0$.

\begin{prop} \label{p:first}
A tight contact structure on $M(-1;r_1,r_2,r_3)$ with matrix of signs
\[
\left ( \begin{matrix} 
1 & 0 & q^3_0 \\
q^1_1 & q^2_1 & q^3_1 
\end{matrix} \right ) 
\]
is isotopic to either a tight contact structure with matrix of signs 
\[
\left ( \begin{matrix} 
0 & 1 & q^3_0 + 1 \\
q^1_1 + 1 & q^2_1  & q^3_1
\end{matrix} \right ) 
\]
or to a tight contact structure with matrix of signs
\[
\left ( \begin{matrix} 
0 & 1 & q^3_0 - 1 \\
q^1_1  & q^2_1-1  & q^3_1
\end{matrix} \right )
\]
provided that the expressions are defined, and all the further basic 
slice decompositions are the same. Here we use the convention 
$\infty \pm 1= \infty$.
\end{prop}

\begin{proof}
We start with the proof of the first isotopy.  If $q^3_0 < a_0^3 -1$
then there is a negative basic slice in the decomposition of $N_0^3$,
therefore we can arrange the basic slice decomposition of $U_3$ so
that the outermost basic slice is negative.  We recall that $M
\setminus (U_1 \cup Z_2 \cup Z_3)$ has boundary slopes $\infty$, $-1$,
and $-1$. Applying Lemma \ref{l:scambio} to $M \setminus (U_1 \cup Z_2 \cup Z_3)$
 we obtain a positive basic slice $B_1$
such that $U_1''=U_1 \cup B_1$ is a tubular neighbourhood of $F_1$ and
$- \partial (M \setminus U_1'')$ has slope $s(- \partial (M \setminus
U_1''))=1$.  Now we divide the proof in two cases.

{\bf Case 1.} If $r_1 = \frac{1}{2}$ then 
\[
A_1^{-1}
\begin{pmatrix} 1\\1\end{pmatrix}
=\begin{pmatrix}
\phantom{-}1 & -1\\
-1 & \phantom{-}2
\end{pmatrix}=
\begin{pmatrix}0\\1
\end{pmatrix},
\]
therefore $s(\del U_1'')=\infty$. It follows that $U_1''$ is the
standard neighbourhood of a destabilization of $F_1$ with twisting
number $0$. Stabilize $F_1$ again, and remove a standard neighbourhood
$U_1'$ of the stabilized curve. We can choose the sign of the
stabilization so that $U_1'' \setminus U_1'$ is a negative basic
slice. The Gluing Theorem \cite[Theorem~4.25]{H1} implies that $U_1'
\setminus V_1'=N_0^1$ is also a negative basic slice, because it glues
to $U_1'' \setminus U_1'$, which is a negative basic slice, to give a
tight contact structure on $U_1'' \setminus V_1'$ with boundary slopes
$-1$ and $\infty$. For this reason $q_0^1$ changes from $q_0^1=1$ to
$q_0^1=0$.  Notice that, in this case, $V_1'=Z_1$. Using Lemma
\ref{l:struttura} applied to $M \setminus (U_1' \cup Z_2 \cup Z_3)$ we
obtain basic slices $B_2'$ and $B_3'$ with boundary slopes $-1$ and
$\infty$. The basic slices $B_2'$ and $B_3'$ are positive by Lemma
\ref{l:backward}, therefore $q_0^2$ changes from $q_0^2=0$ to
$q_0^2=1$, and $q_0^3$ changes to $q_0^3+1$.  This changes the first
row of the matrix of signs from $\left (
\begin{matrix} 1 & 0 & q^3_0 \end{matrix} \right)$ to 
$\left (\begin{matrix} 0 & 1 & q^3_0+1 \end{matrix} \right)$. 

{\bf Case 2.} Suppose now that $r_1 > \frac 12$. Observe that, in view of
Lemma~\ref{l:Zi}, $B_1 \cup N_0^1$ has boundary slopes $0$ and $1$.
Since $B_1$ and $N_0^1$ are both positive basic slices,
by~\cite[Theorem~4.25]{H1} $B_1 \cup N_0^1$ is a positive basic slice
as well. Since 
\[
A_1^{-1}
\begin{pmatrix}1\\1\end{pmatrix}=
\begin{pmatrix}
-\beta_1'-\alpha_1'\\
\beta_1+\alpha_1
\end{pmatrix},
\]
it follows that $B_1$  has boundary slopes 
\[
-\frac{\alpha_1}{\alpha_1'}\quad\text{and}\quad 
-\frac{\beta_1+\alpha_1}{\beta_1'+\alpha_1'} 
\]
in the basis of $-\partial (M \setminus U_1)$.  Moreover,
\[
-\frac{\alpha_1}{\beta_1+\alpha_1}= -\frac{1}{r_1} = 
[a_0^1, \ldots ,a_{k_1}^1] 
\]
implies by Lemma \ref{l:slopeUi} that 
\[
- \frac{\alpha_1}{\alpha_1'}=[a^1_{k_1}, \ldots ,a^1_0],
\]
and by an inductive argument over $k_1$ as in~\cite[Lemma~A4]{OW} 
that 
\[
-\frac{\beta_1+\alpha_1}{\beta_1'+\alpha_1'}= 
[a_{k_1}^1, \ldots, a_1^1]. 
\]
This implies that $B_1\cup N_0^1\cup N^1_1$ is a continued fraction block, 
therefore the signs of the basic slices in $B_1\cup N_0^1\cup N^1_1$
can be shuffled. If we shuffle the sign of the positive basic slice
$B_1 \cup N_0^1$ with the sign of a negative basic slice in $N^1_1$ we
obtain the claimed isotopy.

The proof of the second isotopy is analogous: if $q^3_0>0$, we can
arrange the basic slice decomposition of $U_3$ so that the outermost
basic slice is positive. Applying Lemma \ref{l:scambio} to $M
\setminus (Z_1 \cup U_2 \cup Z_3)$ and proceeding as above we obtain
the second isotopy.
\end{proof} 

\begin{cor}\label{c:ub1}
The number of distinct tight contact structures  on $M(-1;r_1,r_2,r_3)$ 
with $q_0^1 \neq q_0^2$ is bounded above by 
\[
\left(2(a_1^1-1)(a_1^2-1) + (a_0^3-1)(a_1^1+a_1^2-2)\right) (a_1^3-1)
\prod_{i=1}^3 \prod_{j\geq 2} (a_j^i-1)
\]
when $r_2>\frac 12$, by 
\[
\left(2(a_1^1 -1) + (a_0^3 -1)\right) (a_1^3-1)
\prod_{i\neq 2}\prod_{j\geq 2}(a_j^i-1)
\]
when $r_1>r_2=\frac 12$, and by 
\[
2\prod_{j\geq 1}(a_j^3-1)
\]
when $r_1=r_2=\frac 12$. In the above formulae $a_j^i=2$ by convention
if $j>k_i$.
\end{cor}

\begin{proof}
By Proposition \ref{p:first} we can always assume
$(q_0^1,q_0^2)=(1,0)$ unless $(q_0^1,q_0^2)=(0,1)$ and one of the
following cases occurs:
\begin{enumerate}
\item $q_0^3=0$ and $q_1^2=a_1^2-2$,
\item $q_0^3= a_0^3-1$ and $q_1^1= 0$,
\item $q_1^1=0$ and $q_1^2= a_1^2-2$.
\end{enumerate}
There are at most  
\[
(a_1^1-1)(a_1^3-1) \prod_{i=1}^3 \prod_{j\geq 2}(a_j^i-1) 
\]
isotopy classes of tight contact structures in Case (1), 
\[
(a_1^2-1)(a_1^3-1) \prod_{i=1}^3 \prod_{j\geq 2}(a_j^i-1), 
\]
in Case (2), and
\[
a_0^3(a_1^3-1) \prod_{i=1}^3\prod_{j\geq 2}(a_j^i-1) 
\]
in Case (3). But we counted twice the configurations of signs  belonging 
to the groups (1) and (3) or (2) and (3) simultaneously, so 
we have to subtract
\[
2(a_1^3-1) \prod_{i=1}^3 \prod_{j\geq 2}(a_j^i-1)
\] 
from the sum of the above expressions. This shows that the maximum
number of tight contact structures when $(q_0^1,q_0^2)=(0,1)$ and
$r_2>\frac 12$ is
\begin{equation}\label{e:case1}
(a_0^3+a_1^1+a_1^2-4)(a_1^3-1) \prod_{i=1}^3 \prod_{j\geq 2}
(a_j^i-1).
\end{equation}

If $r_1>r_2=\frac 12$, only Case (2) above is possible, so the upper
bound when $(q_0^1,q_0^2)=(0,1)$ is
\begin{equation}\label{e:case1:1}
(a_1^3-1)\prod_{i\neq 2} \prod_{j\geq 2}(a_j^i-1).
\end{equation}
If $r_1=r_2=\frac12$, then none of the above cases can 
occur, and we may always assume $(q_0^1,q_0^2)=(1,0)$. 

Now we consider the case when $(q_0^1,q_0^2)=(1,0)$. By
Proposition~\ref{p:first} the contact structure with matrix of signs
\[
\left ( \begin{matrix} 1 & 0 &
q^3_0 \\ q^1_1 & q^2_1 & q^3_1
\end{matrix} \right ) 
\]
is isotopic to one with matrix 
\[
\left ( \begin{matrix} 1 & 0 & q^3_0 \pm
2\\ q^1_1 \pm 1& q^2_1 \pm 1 & q^3_1
\end{matrix} \right ) 
\]
where the same sign must be chosen in each entry, assuming that
all the expressions are defined. Therefore, we can always assume that
one of the following cases holds:
\begin{enumerate}
\item[(4)] $q_0^3=0$,
\item[(5)] $q_0^3=1$,
\item[(6)] $q_0^3 \neq 0,1$ and $q_1^1=0$,
\item[(7)] $q_0^3 \neq 0, 1$ and $q_1^2=0$.
\end{enumerate} 
Each one of Cases (4) and (5) allows the existence of at most 
\[
(a_1^1-1)(a_1^2-1)(a_1^3-1) \prod_{i=1}^3 \prod_{j\geq 2}
(a_j^i-1)
\]
distinct isotopy classes of contact structures, Case (6) allows 
\[
(a_0^3-2)(a_1^2-1)(a_1^3-1) \prod_{i=1}^3 \prod_{j\geq 2}(a_j^i-1), 
\]
and Case (7) allows
\[
(a_0^3-2)(a_1^1-1)(a_1^3-1) \prod_{i=1}^3 \prod_{j\geq 2}(a_j^i-1). 
\]
However, we have counted twice the contact structures with
$(q_1^1,q_1^2)=(0,0)$ belonging to both Cases (6) and (7). 
Therefore, we must subtract
\[
(a_0^3-2)(a_1^3-1) \prod_{i=1}^3 \prod_{j\geq 2}(a_j^i-1). 
\]
Thus, when $r_2>\frac 12$ the total number of potential tight contact
structures with $(q_0^1,q_0^2)=(1,0)$ is
\begin{equation}\label{e:case2}
\left( 2(a_1^1 -1)(a_1^2 -1) + (a_0^3 -2)(a_1^1 + a_1^2 -3)\right) (a_1^3-1)
\prod_{i=1}^3 \prod_{j\geq 2}(a_j^i-1). 
\end{equation}
Adding up~\eqref{e:case1} and~\eqref{e:case2} we obtain the stated
formula in the case $r_2>\frac 12$.

If $r_1>r_2=\frac 12$ then  Case (7) cannot
occur because $q_1^2=\infty$, therefore the total number of potential tight 
contact structures with $(q_0^1,q_0^2)=(1,0)$ is
\begin{equation}\label{e:case2:1}
\left(2(a_1^1 -1) + (a_0^3 -2)\right) (a_1^3-1)
\prod_{i\neq 2} \prod_{j\geq 2}(a_j^i-1). 
\end{equation}
Adding up~\eqref{e:case1:1} and~\eqref{e:case2:1} gives the statement
in this case.

When $r_1=r_2=\frac 12$ only Cases (4) and (5) can occur, giving the
upper bound
\[
2 \prod_{j\geq 1}(a_j^3-1),
\]
which coincides with the stated formula. 
\end{proof}

Next we turn to the second possibility, when $q^1_0=q^2_0$.

\begin{prop}\label{p:rest}
Let $\xi$ is a tight contact structure on $M=M(-1; r_1, r_2,r_3)$ such that 
$q_0^1=q_0^2$.  If $(q^1_0,q^2_0)=(1,1)$ then $q^3_0=0$ and if $(q^1_0,q^2_0)=(0,0)$ then 
$q^3_0= a^3_0 -1$.
\end{prop}

\begin{proof}
Recall that $M \setminus (Z_1 \cup Z_2 \cup  U_3)$ has boundary slopes
$0$, $-1$, and $\infty$. Applying Lemma~\ref{l:scambio} to $M
\setminus (Z_1 \cup Z_2 \cup U_3)$  we get a basic
slice $B_3$ with the sign opposite to the sign of $N_0^1$ and $N_0^2$.
Then $U_3''=B_3 \cup U_3$ is a tubular neighbourhood of $F_3$ so that
$- \partial (M \setminus U_3'')$ has slope $0$. Since
\[
A_3^{-1}\begin{pmatrix} 1\\0\end{pmatrix}=
\begin{pmatrix}-\be'_3\\\be_3\end{pmatrix}, 
\]
$U_3''$ has boundary slope $-\frac{\be_3}{\be'_3}$.  Now suppose that
$r_3\neq\frac1{a_0^3}$. Then, by induction on $k_3$ as
in~\cite[Lemma~A4]{OW}, we have
\[
-\frac{\beta_3}{\beta_3'}=[a_{k_3}^3, \ldots ,a_1^3].
\]
Since in the basis of $- \partial (M \setminus V_3)$ the toric
annulus $B_3 \cup N_0^3$ has boundary slopes 
\[
[a^3_{k_3}, \ldots, a^3_1-1]\quad\text{and}\quad 
[a^3_{k_3}, \ldots ,a^3_1], 
\]
which are joined by an
edge in the Farey Tessellation, by~\cite[Theorem~4.25]{H1} $B_3 \cup
N_0^3$ is a basic slice, and it is tight if and only if $B_3$ and all
the basic slices in the basic slice decomposition of $N_0^3$ have the
same sign.  This happens if and only if
\[
(q^1_0,q^2_0,q_0^3)=(1,1,0)\quad\text{or}\quad
(q^1_0,q^2_0,q_0^3)=(0,0,a^3_0 -1).
\]
Now suppose that  $r_3 = \frac{1}{a_0^3}$. Since
\[
A_3^{-1}\begin{pmatrix} 1\\0\end{pmatrix} = 
\begin{pmatrix}
a_0^3 & a_0^3-1\\
1 & \phantom{-1}1
\end{pmatrix}^{-1}
\begin{pmatrix} 1\\0\end{pmatrix} = 
\begin{pmatrix} 1\\1\end{pmatrix},
\]
$\partial U_3''$ has slope $1$. Therefore, in this case $B_3 \cup N_0^3$ has boundary
slopes $-1$ (as observed after Lemma~\ref{l:slopeUi}) and $1$. Thus, we 
can argue as in the case $r_3\neq\frac1{a_0^3}$ and draw the same conclusion.  
\end{proof}

\begin{prop} \label{p:second}
Two tight contact structures on $M(-1;r_1,r_2,r_3)$ with matrices of signs
\[
\left ( \begin{matrix} 
1 & 1 & 0 \\
 q^1_1 &  q^2_1 &  q^3_1
\end{matrix} \right )
\]
and
\[
\left ( \begin{matrix} 
0 & 0 &  a^3_0 -1 \\
q^1_1 &  q^2_1 & q^3_1 + 1
\end{matrix} \right )
\]
and with identical basic slice decomposition in all the further 
continued fraction blocks are isotopic, whenever the symbols are
 defined. Here we use the convention $\infty\pm 1=\infty$.
\end{prop}

\begin{proof}
Suppose that $q^1_0 = q^2_0$.  Applying Lemma \ref{l:scambio} to $M
\setminus (Z_1 \cup Z_2 \cup U_3)$ we obtain a basic slice $B_3$ such
that $U_3''= B_3 \cup U_3$ is a tubular neighbourhood of $F_3$, and $-
\partial (M \setminus U_3'')$ has slope $s(- \partial (M \setminus
U_3''))=0$. Now we divide the proof into two cases.

{\bf Case 1.} Suppose first that $r_3 = \frac{1}{a_0^3}$. As in the
proof of Proposition~\ref{p:rest}, we see that $U_3''$ has slope $1$,
and therefore it is the standard neighbourhood of a destabilization of
$F_3$ with twisting number $1$. Stabilize $F_3$ again, and remove a
standard neighbourhood $U_3'$ of the stabilized curve. We can choose
the sign of the stabilization so that $U_3'' \setminus U_3'$ is a
negative basic slice.  Proceeding as in Case~1 of the proof of
Proposition \ref{p:first}, we can change the first row of the matrix
of signs from $\left (
\begin{matrix} 1 & 1 & 0 \end{matrix} \right)$ to $\left (
\begin{matrix} 0 & 0 & a_0^3-1 \end{matrix} \right)$.

{\bf Case 2.} If $r_3 \neq \frac{1}{a_0^3}$ then $B_3 \cup N_0^3$ is a
basic slice with boundary slopes $[a^3_{k_3}, \ldots ,a^3_1-1]$ and
$[a^3_{k_3}, \ldots ,a^3_1]$ , and $B_3$ and all the basic slices in
$N_0^3$ have the same sign (cf.~the proof of
Proposition~\ref{p:rest}).  $B_3 \cup N_0^3 \cup N_1^3$ has boundary
slopes $[a^3_{k_3}, \ldots ,a_2^3-1]$ and $- \frac{\beta_3}{\beta_3'}=
[a^3_{k_3}, \ldots ,a_1^3]$ computed in the basis of $- \partial (M
\setminus V_3)$.  This implies that $B_3 \cup N_0^3 \cup N^3_1$ is a
continued fraction block. According to~\cite[Lemma~4.14]{H1} we can
swap the signs of the basic slice $B_3\cup N_0^3$ and the sign of a
basic slice of $N^3_1$.  This gives the stated isotopy.
\end{proof}

\begin{cor}\label{c:ub2}
The number of isotopy classes of tight contact structures with $q_0^1=q_0^2$ 
carried by $M(-1;r_1,r_2,r_3)$ is bounded above by
\[
(a_1^1-1)(a_1^2-1)a_1^3 \prod _{i=1}^3\prod_{j\geq 2}(a_j^i -1)
\]
if $r_3\neq\frac1{a_0^3}$, and by
\[
\prod _{i=1}^2\prod _{j\geq 1}(a_j^i -1)
\]
if $r_3 = \frac1{a_0^3}$. In the above formulae $a_j^i=2$ by convention
if $j>k_i$.
\end{cor}

\begin{proof}
By Proposition~\ref{p:rest} there are two possibilities for the first
row of the matrix of signs $(q_0^1, q_0^2, q_0^3)$ defined by the
number of positive basic slices in the outermost continued fraction
blocks $N_0^i$: either $(q_0^1, q_0^2, q_0^3)=(1,1,0)$ or $(q_0^1,
q_0^2, q_0^3)=(0,0,a_0^3 -1)$.  By Proposition~\ref{p:second} each
potentially tight contact structure with $(0,0,a_0^3 -1)$ as first row
of the matrix of signs is isotopic to one with $(1,1,0)$ unless
$q_1^3=0$. If $r_3=\frac1{a_0^3}$ we have $q_1^3=\infty$, therefore 
in this case we get the bound
\[
(a_1^1-1)(a_1^2-1)
\prod _{i=1}^2\prod _{j\geq 2}(a_j^i -1) = 
\prod _{i=1}^2\prod _{j\geq 1}(a_j^i -1).
\]
When $r_3\neq\frac1{a_0^3}$ we can consider two
cases:
\begin{enumerate}
\item $(q_0^1,q_0^2,q_0^3)=(1,1,0)$, or
\item $(q_0^1,q_0^2,q_0^3)=(0,0,a_0^3-1)$ and $q_1^3=0$. 
\end{enumerate}
Case (1) gives the upper bound
\begin{equation}\label{e:first}
(a_1^1-1)(a_1^2-1)(a_1^3-1) \prod _{i=1}^3\prod_{j\geq 2}(a_j^i -1) 
\end{equation}
and Case (2) gives
\begin{equation}\label{e:second}
(a_1^1-1)(a_1^2-1) \prod_{i=1}^3\prod_{j\geq 2}(a_j^i -1).
\end{equation}
Adding up~\eqref{e:first} and~\eqref{e:second} gives the 
statement when $r_3\neq\frac1{a_0^3}$. 
\end{proof}

\section{Contact Ozsv\'ath--Szab\'o invariants}\label{s:OSz}
\sh{Ozsv\'ath--Szab\'o homology groups}

In a remarkable series of papers \cite{OSzF1, OSzF2, OSzF4, OSz6}
Ozsv\'ath and Szab\'o defined new invariants of many low--dimensional
objects --- including contact structures on closed $3$--manifolds.
Heegaard Floer theory associates to a closed, oriented spin$^c$ 
$3$--manifold
$(Y, \t)$ the abelian groups $\hf (Y, \t)$, $HF^\infty(Y,\t)$,
$HF^-(Y,\t)$ and $HF^+(Y,\t)$, called the \emph{Ozsv\'ath--Szab\'o
homology groups}. If $(W, \s)$ is an oriented spin$^c$ cobordism
between two spin$^c$ $3$--manifolds $(Y_1, \t_1)$ and $(Y_2, \t_2)$
and $HF^\bullet (Y,\t_i)$, $i=1,2$ is any of the groups above, there is a
homomorphism
\[
 F^\bullet_{W, \s}\co HF^\bullet(Y_1, \t_1)\to HF^\bullet(Y_2, \t_2).
\]
In this paper we shall be mainly concerned with the groups
$\hf(Y,\t)$, which are always finitely generated. The symbol $\hf (Y)$
will denote the direct sum
\[
\hf(Y):=\oplus_{\t}\hf (Y,\t)
\]
over all spin$^c$ structures $\t$ on $Y$. Since there are only
finitely many spin$^c$ structures with nonvanishing $\hf$--group, $\hf
(Y)$ is still finitely generated. A rational homology 3--sphere
$Y$ is called an \emph{$L$--space} if $\hf (Y, \t)\cong \Z$ for all
$\t \in Spin^c(Y)$.

Let $Y$ be a closed, oriented 3--manifold and let $K\subset Y$ be a
framed knot with framing $f$. Let $Y(K)$ denote the 3--manifold given
by surgery along $K\subset Y$ with respect to the framing $f$ and
$Y'(K)$ the 3-manifold we get by performing surgery along $K$ with
framing $f+\mu$, where $\mu$ denotes the meridian of $K$. The
surgeries determine cobordisms $X_1$ from $Y$ to $Y(K)$, $X_2$ from
$Y(K)$ to $Y'(K)$ and $X_3$ from $Y'(K)$ back to $Y$. The following
result can be deduced (cf.~the discussion at the beginning
of~\cite[Section~3]{OSzabs} and~\cite[page 934]{LSuj})
from~\cite[Theorem~9.16]{OSzF2} and~\cite[Subsection~4.1]{OSzF4}.

\begin{thm}[Surgery exact triangle]
\label{t:triangle}
The Ozsv\'ath--Szab\'o homology groups of $Y$, $Y(K)$ and $Y'(K)$ fit
into an exact triangle
\[
\begin{graph}(6,2)
\graphlinecolour{1}\grapharrowtype{2}
\textnode {A}(1,1.5){$\hf (Y)$}
\textnode {B}(5, 1.5){$\hf (Y(K))$}
\textnode {C}(3, 0){$\hf (Y'(K))$}
\diredge {A}{B}[\graphlinecolour{0}]
\diredge {B}{C}[\graphlinecolour{0}]
\diredge {C}{A}[\graphlinecolour{0}]
\freetext (3,1.8){$F_{X_1}$}
\freetext (4.6,0.6){$F_{X_2}$}
\freetext (1.4,0.6){$F_{X_3}$}
\end{graph}
\]
where 
\[
F_{X_i} = \sum_{\s\in Spin^c (X_i)} \pm F_{X_i,\s},\quad i=1,2,3.
\]
\qed
\end{thm}

It was proved in~\cite{OSzF1, OSzabs} that for each spin$^c$ structure
$\t$ the Ozsv\'ath--Szab\'o homology group $\hf(Y,\t)$ comes with a
natural relative $\Z/div(\t)\Z$--grading, where $div(\t)$ is the
divisibility of $c_1(\t)$ in $H^2(Y;\Z)$. If $\t \in Spin ^c(Y)$ is
torsion, that is, $c_1(\t )\in H^2 (Y; \Z )$ is a torsion element,
then $div(\t)=0$, and therefore $\hf(Y,\t)$ has a natural relative
$\Z$--grading. This relative $\Z$--grading admits a natural lift to an
absolute $\Q$--grading.  In conclusion, for a torsion spin$^c$
structure $\t$ the Ozsv\'ath--Szab\'o homology group $\hf (Y, \t)$
splits as
\[
\hf (Y, \t )=\oplus _{n\in\Z}\hf_{d_0+n}(Y,\t), 
\]
where the degree $d_0\in\Q$ is determined mod 1 by $\t$.  Moreover,
when $\t \in\Spin^c(Y)$ has torsion first Chern class, there is an
isomorphism between the homology groups $\hf_d (Y,\t)$ and $\hf
_{-d}(-Y,\t)$.

Let $(W, \s )$ be a spin$^c$ cobordism between two spin$^c$ manifolds
$(Y_1,\t_1) $ and $(Y_2, \t _2)$.  If the spin$^c$ structures $\t_i$
are both torsion and $x\in\hf (Y_1, \t_1)$ is a homogeneous element of
degree $d(x)$, then $F_{W,\s}(x)\in \hf (Y_2, \t _2)$ is also
homogeneous of degree
\begin{equation}\label{e:d-shift}
d(x) + \frac{1}{4}(c_1^2(\s )-3\sigma (W)-2\chi (W)).
\end{equation}
We need one more piece of information. Recall that the set of spin$^c$
structures comes equipped with a natural involution, usually denoted
by $\t \mapsto \ot$. The spin$^c$ structure $\ot$, called
the~\emph{conjugate} of $\t$, is defined as follows: If one thinks of
a spin$^c$ structure as a suitable equivalence class of nowhere zero
vector fields (cf. \cite{OSzF1}) then the above involution is the map
induced by multiplying a representative vector field by
$(-1)$. 

\begin{thm}[\cite{OSzF2}, Theorem~2.4]\label{t:conj-iso}
There is a natural isomorphism 
\[
\J_Y\co\hf (Y,\t)\to\hf (Y,\ot)
\]
\qed
\end{thm}

A spin$^c$ structure $\t\in Spin^c(Y)$ is induced by a \emph{spin}
structure exactly when $c_1(\t)=0$, or equivalently when $\t = \ot$.
According to~\cite[Theorem~3.6]{OSzF4}, given a spin$^c$ cobordism
$(W,\s)$ we have
\begin{equation}\label{e:szam}
\J_{Y'}\circ F_{W,\s} = F_{W,\os}\circ\J_Y,
\end{equation}
where $\os$ is the spin$^c$ structure on the 4--manifold $W$ conjugate
to $\s$. This means that, if one thinks of $\s\in Spin^c (W)$ as a
suitable equivalence class of almost--complex structures defined on
$W\setminus\{\text{finitely many points}\}$, then $\os$ is represented
by the negative $-J$ of any almost--complex structure $J$ representing
$\s$.

\sh{Contact $(\pm 1)$--surgery}

Suppose that $L\subset (Y, \xi )$ is a Legendrian knot in a contact
3--manifold. Let $Y_L^{\pm }$ denote the 3--manifold obtained by doing
$(\pm 1)$--surgery along $L$, where the surgery coefficient is
measured with respect to the contact framing of $L$. According to the
classification of tight contact structures on a solid torus~\cite{H1},
the contact structure $\xi\vert _{Y-\nu L}$ extends uniquely (up to
isotopy) to the surgered manifolds $Y_L^+$ and $Y_L^-$ as a tight
structure on the glued--up torus. Therefore, the knot $L$ with a
$(+1)$ or $(-1)$ on it uniquely specifies a contact 3--manifold
$(Y_L^+, \xi _L^+)$ or $(Y_L^-, \xi _L^-)$.  (For more about contact
surgery see \cite{DG1, DG2, DGS}.)  In particular, a Legendrian link
${\mathbb {L}} \subset (S^3, \xi _{st})$ in the standard contact
3--sphere (which can be represented by its front projection) defines a
contact structure once the surgery coefficients $(+1)$ and $(-1)$ are
specified on its components.

\sh{Contact Ozsv\'ath--Szab\'o invariants}

In~\cite{OSz6} Ozsv\'ath and Szab\'o define an invariant
\[
c(Y, \xi)\in \hf (-Y,\t_{\xi })/\{\pm 1\}
\]
assigned to a positive, cooriented contact structure $\xi$ on $Y$.  In
fact, $\xi$ (as an oriented 2--plane field) determines an element
$(d(\xi),\t _{\xi })\in \Ht$ and according to \cite{OSz6} the contact
invariant $c(Y,\xi )$ is an element of $\hf
_{-d(\xi)}(-Y,\t_{\xi})/\{\pm 1\}$. Moreover, if $c_1(\xi )\in H^2 (Y;
\Z )$ is torsion then
\begin{equation}\label{e:d3-inv}
d(\xi)= \frac{1}{4}(c_1^2(X, J)- 3\sigma (X) -2 \chi (X)+2),
\end{equation}
where $(X,J)$ is a compact almost--complex 4--manifold with $\partial
X=Y$, and $\xi $ is homotopic to the distribution of complex
tangencies on $\partial X$.

The main properties of the contact Ozsv\'ath--Szab\'o invariant are
summarized in the following two theorems.

\begin{thm}[\cite{OSz6}] \label{t:item1}
If $(Y, \xi )$ is overtwisted, then $c(Y, \xi )=0$. If $(Y, \xi )$ is
Stein fillable then $c(Y, \xi )\neq 0$.  In particular, for the
standard contact structure $(S^3, \xi _{st})$ the invariant $c(S^3,
\xi _{st})$ is nonzero. \qed
\end{thm}

Given a spin$^c$ cobordism $(W,\s)$ between spin$^c$ 3--manifolds
$(Y_1,\t_1)$ and $(Y_2,\t_2)$, the homomorphism $F_{W,\s}$ clearly
induces a map between the sets $\hf(Y_i,\t_i)/\{\pm 1\}$,
$i=1,2$. Likewise, for any spin$^c$ 3--manifold $(Y,\t)$ the isomorphism 
$\J_Y\co\hf(Y,\t)\to\hf(Y,\t)$ induces a map
\[
\hf(Y,\t)/\{\pm\}\to\hf(Y,\t)/\{\pm\}.
\]
Abusing notation, throughout the paper we shall keep denoting
such maps by $F_{W,\s}$ and $\J_Y$.

\begin{thm}[\cite{LS3, OSz6}]\label{t:item2}
Suppose that $(Y_2,\xi_2)$ is obtained from $(Y_1,\xi _1)$ by a
contact $(+1)$--surgery, and let $-W$ be the cobordism induced by the
surgery with reversed orientation. Then, 
\[
F_{-W} (c(Y_1, \xi_1))= c(Y_2,\xi_2). 
\]
In particular, if $c(Y_2, \xi_2)\neq 0$ then $(Y_1, \xi_1)$ is
tight. \qed
\end{thm}

Since by~\cite[Proposition~8]{DG1} contact $(+1)$--surgery along a
Legendrian knot $L$ is cancelled by a contact $(-1)$--surgery along a
Legendrian push--off of $L$, Theorem~\ref{t:item2} immediately
implies:

\begin{cor}\label{c:legsurg}
If $(Y_2,\xi _2)$ is obtained by Legendrian surgery along a Legendrian
knot in $(Y_1, \xi _1)$ and $c(Y_1, \xi _1)\neq 0$, then $c(Y_2, \xi
_2)\neq 0$. In particular, $(Y_2, \xi _2)$ is tight. \qed
\end{cor}

\noindent An easy application of the surgery exact triangle together
with Theorem~\ref{t:item2} gives

\begin{lem}[\cite{LS3}, Lemma~2.5] \label{l:etak}
The contact structure $\eta _1$ on $S^1\times S^2$ given as
$(+1)$--surgery on a Legendrian unknot with Thurston--Bennequin number
$-1$ in $(S^3,\xi_{st})$ has nonvanishing contact Ozsv\'ath--Szab\'o
invariant. \qed
\end{lem}

\section{Tight contact structures on $M(-1;\frac12,\frac12,\frac1p)$} 
\label{s:special}

In this section we define three contact structures $\xi_1$, $\xi_2$
and $\Xi$ on the 3--manifold $M_p=M(-1;\frac{1}{2}, \frac{1}{2},
\frac{1}{p})$ for each $p\geq 2$, we prove that they are distinct up
to homotopy, that $\xi_1$ and $\xi_2$ are Stein fillable and that
$\Xi$ has nonzero Ozsv\'ath--Szab\'o contact invariant. Combined with
the results of Section~\ref{s:convex}, this gives the complete
classification of tight contact structures on $M_p$ for every $p$.  In
the last subsection we show that for $p>2$ the contact 3--manifold
$(M_p,\Xi)$ is not Stein fillable, and for $p\not\equiv 2\bmod 8$ is
not symplectically fillable.

\sh{Heegaard Floer groups of $M_p$}

The oriented manifold $-M_p$ is represented by the third surgery
diagram of Figure~\ref{f:kirby}. The three dotted circles denoted by
$\mu_a$, $\mu_b$ and $\mu_c$ represent (up to sign) elements of
$H_1(-M_p;\Z)$. It is easy to check that
\[
H_1(M_p;\Z)=
\begin{cases}
\langle \mu_b\ |\ 4\mu_b=0\rangle\cong \Z/4\Z\quad\text{for $p$ odd},\\
\langle \mu_b,\mu_c\ |\ 2\mu_b=2\mu_c=0\rangle\cong\Z/2\Z\oplus\Z/2\Z
\quad\text{for $p$ even},
\end{cases}
\]
In particular, $M_p$ has four spin$^c$ structures for every $p$.

\begin{figure}[ht]
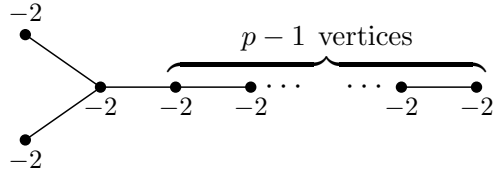

\begin{center}
\begin{graph}(10,2)(0,-1)
\graphnodesize{0.15}

  \roundnode{top}(1,0.7)
  \roundnode{bottom}(1,-0.7)
  \roundnode{n0}(2,0)
  \roundnode{n11}(3,0)
  \roundnode{n12}(4,0)

  \roundnode{n13}(6,0)
  \roundnode{n14}(7,0)

  \freetext(5,0.30){$\overbrace{\hspace{120pt}}$}
  \freetext(5,0.65){$p-1$ vertices}

  \edge{top}{n0}
  \edge{bottom}{n0}
  \edge{n0}{n11}
  \edge{n11}{n12}
  \autonodetext{n12}[e]{\large $\cdots$}

  \autonodetext{n13}[w]{\large $\cdots$}
  \edge{n13}{n14}

  \autonodetext{top}[n]{\small $-2$}
  \autonodetext{bottom}[s]{\small $-2$}
  \autonodetext{n0}[s]{\small $-2$}
  \autonodetext{n11}[s]{\small $-2$}
  \autonodetext{n12}[s]{\small $-2$}
  \autonodetext{n13}[s]{\small $-2$}
  \autonodetext{n14}[s]{\small $-2$}

\end{graph}
\end{center}
\caption{The 4--dimensional plumbing $P$ with boundary $-M_p$}
\label{f:D} 
\end{figure}

The sequence of Kirby calculus moves going from the third to the
seventh diagram of Figure~\ref{f:kirby} shows that $-M_p$ is the
boundary of $P$, the plumbing of spheres given by Figure~\ref{f:D}.
This amounts to saying that $-M_p$ is the link of the singularity
$D_{p-2}$.

\begin{figure}[ht]
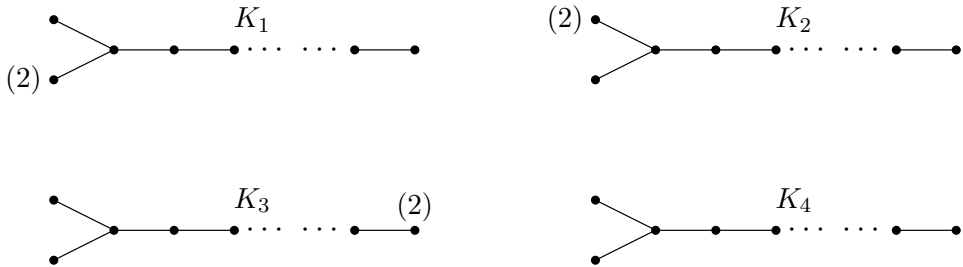

\unitlength=.8cm
\begin{center}
\begin{graph}(18,4.5)(0,0)
\graphnodesize{0.15}

  \roundnode{top1}(1,1)
  \roundnode{bottom1}(1,0)
  \roundnode{n1}(2,0.5)
  \roundnode{n2}(3,0.5)
  \roundnode{n3}(4,0.5)

  \roundnode{n4}(6,0.5)
  \roundnode{n5}(7,0.5)

  \freetext(4.3,1){$K_3$}

  \edge{top1}{n1}
  \edge{bottom1}{n1}
  \edge{n1}{n2}
  \edge{n2}{n3}
  \autonodetext{n3}[e]{\large $\cdots$}

  \autonodetext{n4}[w]{\large $\cdots$}
  \edge{n4}{n5}
  \autonodetext{n5}[n]{$(2)$}

  \roundnode{top2}(10,1)
  \roundnode{bottom2}(10,0)
  \roundnode{m1}(11,0.5)
  \roundnode{m2}(12,0.5)
  \roundnode{m3}(13,0.5)

  \roundnode{m4}(15,0.5)
  \roundnode{m5}(16,0.5)

  \freetext(13.3,1){$K_4$}

  \edge{top2}{m1}
  \edge{bottom2}{m1}
  \edge{m1}{m2}
  \edge{m2}{m3}
  \autonodetext{m3}[e]{\large $\cdots$}

  \autonodetext{m4}[w]{\large $\cdots$}
  \edge{m4}{m5}

  \roundnode{top3}(1,4)
  \roundnode{bottom3}(1,3)
  \roundnode{l1}(2,3.5)
  \roundnode{l2}(3,3.5)
  \roundnode{l3}(4,3.5)

  \roundnode{l4}(6,3.5)
  \roundnode{l5}(7,3.5)

  \freetext(4.3,4){$K_1$}

  \edge{top3}{l1}
  \edge{bottom3}{l1}
  \autonodetext{bottom3}[w]{$(2)$}
  \edge{l1}{l2}
  \edge{l2}{l3}
  \autonodetext{l3}[e]{\large $\cdots$}

  \autonodetext{l4}[w]{\large $\cdots$}
  \edge{l4}{l5}

  \roundnode{top4}(10,4)
  \roundnode{bottom4}(10,3)
  \roundnode{p1}(11,3.5)
  \roundnode{p2}(12,3.5)
  \roundnode{p3}(13,3.5)

  \roundnode{p4}(15,3.5)
  \roundnode{p5}(16,3.5)

  \freetext(13.3,4){$K_2$}

  \edge{top4}{p1}
  \edge{bottom4}{p1}
  \autonodetext{top4}[w]{$(2)$}
  \edge{p1}{p2}
  \edge{p2}{p3}
  \autonodetext{p3}[e]{\large $\cdots$}

  \autonodetext{p4}[w]{\large $\cdots$}
  \edge{p4}{p5}

\end{graph}
\end{center}
\caption{Initial vectors}
\label{f:vectors} 
\end{figure}

Consider the four $2$--cohomology classes $K_i$, $i=1,\ldots,4$, on
 $P$ whose values on the standard homology generators are given by
 Figure~\ref{f:vectors} (each number in parentheses indicates the
 value of $K_i$ on the homology generator corresponding to the nearby
 vertex, and by convention no number is present if such value is
 zero).

\begin{defn}
Define $\t_i$, for $i=1,\ldots,4$, to be the spin$^c$ structure on
$-M_p$ which is the restriction of the spin$^c$ structure on $P$
specified by the characteristic element $K_i$.
\end{defn}

An easy calculation shows that 
\[
\t_1=\t_4+\mu_b,\quad\t_2=\t_4+\mu_a+\mu_b=\t_4+\mu_c
\quad\text{and}\quad\t_3=\t_4+\mu_a,
\]
where $\mu_a$, $\mu_b$ and $\mu_c$ are the homology classes defined in
Figure~\ref{f:kirby}. This implies that $\{\t_1,\t_2,\t_3,\t_4\}$ is
the whole set of spin$^c$ structures on $-M_p$. When $p$ is even, $M_p$
has $4$ spin structures, therefore in this case each of the spin$^c$
structures $\t_i$ is induced by a spin structure. When $p$ is odd,
$M_p$ carries $2$ spin structures, which induce $\t_3$ and $\t_4$. In
fact, observe that $\t_4$ is always induced by a spin structure
because $K_4=0$ and therefore $c_1(\t_4)=0$. On the other hand, we
always have $\mu_a=\mu_b-\mu_c$, while when $p$ is odd $\mu_c=3\mu_b$
and therefore
\[
c_1(\t _3)=c_1(\t _4)+2\mu_a=-4\mu_b=0.
\]

\begin{prop}\label{p:mp}
We have
\[
\hf (-M_p,\t_1)\cong\hf (-M_p,\t_2)\cong\Z_{(0)},\ \ 
\hf (-M_p,\t_3)\cong\Z_{(\frac{p-2}{4})},\ \ 
\hf (-M_p,\t_4)\cong\Z_{(\frac{p+2}{4})}.
\]
\end{prop}

\begin{proof}
Since $-M_p$ is the boundary of the plumbing $P$, we can apply the
algorithm of~\cite{plum} to determine its Ozsv\'ath--Szab\'o homology
groups. It is easy to check that the four cohomology classes $K_i$,
$i=1,\ldots,4$, of Figure~\ref{f:vectors} provide initial
characteristic vectors in the sense of~\cite{plum}.
An easy computation shows that
\[
K_1^2=K_2^2=-p-2,\quad K_3^2=-4\quad\text{and}\quad K_4^2=0.
\]
The 3--manifold $-M_p$ has elliptic geometry and therefore it is an
$L$--space by~\cite[Proposition~2.3]{OSzlens}. Since by~\cite{plum} 
\[
d(-M_p,\t_i)=\frac{K_i^2+p+2}4,\quad i=1,\ldots,4,
\]
this immediately implies the statement.
\end{proof}

\sh{Contact structures on $M_p$}

\begin{defn}
Let $\xi_1$ and $\xi_2$ be the contact structures defined respectively
by the contact surgery diagrams of Figure~\ref{f:structure}(a)
and~\ref{f:structure}(b).
\begin{figure}[ht]
\begin{center}
\psfrag{A}{\small $+1$}
\psfrag{B}{\small $-1$}
\psfrag{L}{\small $L$}
\psfrag{(a)}{\small $(a)$}
\psfrag{(b)}{\small $(b)$}
\psfrag{p-1}{\small $p-1$}
\psfrag{cusps}{\small left cusps}
\psfrag{p-2}{\small $p-2$}
\includegraphics[height=5cm]{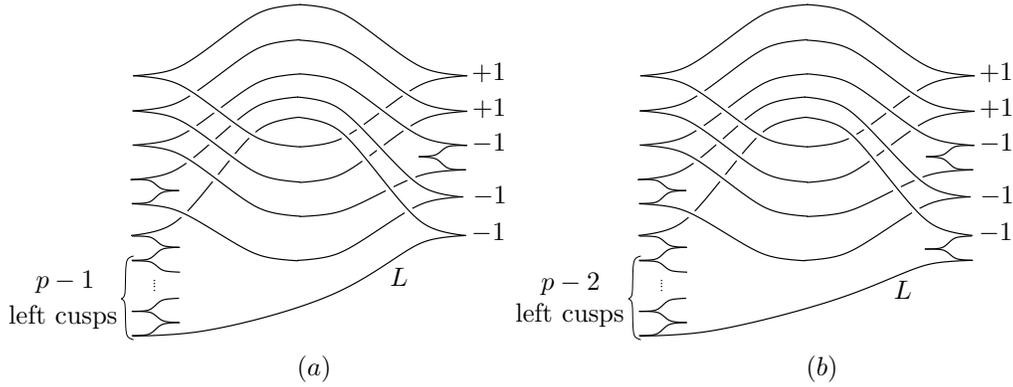}
\end{center}
\caption{\quad The contact structures $\xi_1$ and $\xi_2$ on $M_p$}
\label{f:structure}
\end{figure}
\end{defn}

\begin{defn}
Let $\Xi$ be the contact structure defined by the contact surgery
diagram of Figure~\ref{f:bigksi}.
\begin{figure}[ht]
\begin{center}
\psfrag{A}{\small $+1$}
\psfrag{B}{\small $-1$}
\psfrag{L}{\small $L$}
\psfrag{p-1}{\small $p-1$}
\psfrag{cusps}{\small left cusps}
\includegraphics[height=5cm]{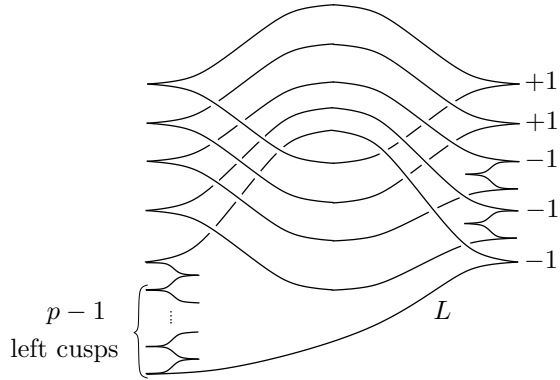}
\end{center}
\caption{The contact structure $\Xi$ on $M_p$}
\label{f:bigksi}
\end{figure}
\end{defn}

\begin{lem}\label{l:spin-c}
We have
\[
\{\t_{\xi_1}, \t_{\xi_2}\}=\{\t_1, \t_2\}\quad\text{and}
\quad 
\t_{\Xi}=\t_3.
\]
\end{lem}

\begin{proof}
The contact surgery presentation of each contact structure
$\xi\in\{\xi_1,\xi_2,\Xi\}$ can be interpreted as a simply connected
4--manifold with boundary, endowed with a characteristic 2--cohomology
class $K_\xi$. The class $K_\xi$ is uniquely determined by requiring
that it evaluates on a 2--homology generator corresponding to a given
Legendrian knot ${\mathcal L}$ in the surgery presentation as the
rotation number of ${\mathcal L}$ (once an orientation for ${\mathcal
L}$ is chosen). Moreover, by~\cite{DGS} the spin$^c$ structure
determined by $K_\xi$ restricts to the spin$^c$ structure associated
to $\xi$ on the boundary.

On the other hand, the given surgery presentation can be viewed
smoothly as the first diagram of Figure~\ref{f:kirby}.  
\begin{figure}
\begin{center}
\psfrag{0}{{\scriptsize $0$}}
\psfrag{3}{\scriptsize $-3$}
\psfrag{-p-1}{\small $-p-1$}
\psfrag{m1}{\scriptsize $-1$}
\psfrag{1}{\scriptsize $1$}
\psfrag{-1}{\scriptsize $-1$}
\psfrag{-2}{\scriptsize $-2$}
\psfrag{2}{\scriptsize $2$}
\psfrag{p-1}{\small $p-1$}
\psfrag{p}{\small $p$}
\psfrag{-p}{\small $-p$}
\psfrag{reverse}{\scriptsize orientation reversal}
\psfrag{mua}{\small $\mu_a$}
\psfrag{mub}{\small $\mu_b$}
\psfrag{muc}{\small $\mu_c$}
\includegraphics[height=0.94\textheight]{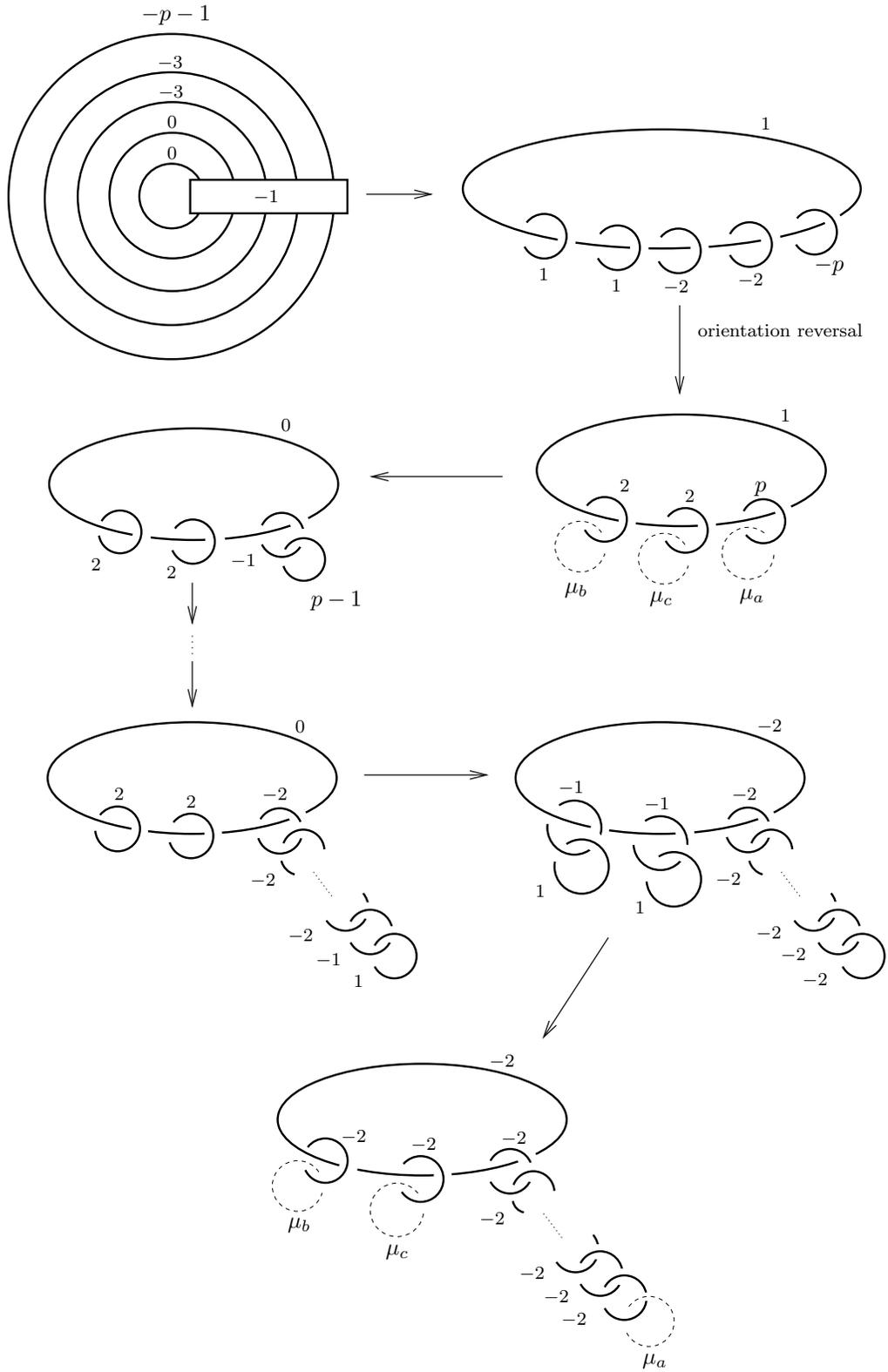}
\end{center}
\caption{Kirby diagrams for $\pm M_p$}
\label{f:kirby}
\end{figure}
By carrying along the class $K_\xi$ during the Kirby moves of
Figure~\ref{f:kirby} (and observing that blowups and blowdowns do not
change the spin$^c$ structure on the boundary) one can check that the
spin$^c$ structures $\t_{\xi_1}$, $\t_{\xi_2}$ and $\t_\Xi$ are the
restrictions to the boundary of the spin$^c$ structures on the
4--dimensional plumbing determined by, respectively, the
characteristic classes $C_1$, $C_2$ and $C_3$ given in
Figure~\ref{f:classes}.
\begin{figure}[ht]
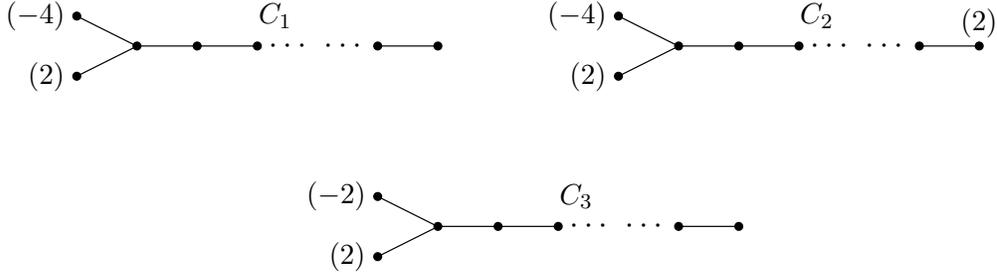

\unitlength=.8cm
\begin{center}
\begin{graph}(18,4.5)(0,0)
\graphnodesize{0.15}
  \roundnode{top1}(6,1)
  \roundnode{bottom1}(6,0)
  \roundnode{n1}(7,0.5)
  \roundnode{n2}(8,0.5)
  \roundnode{n3}(9,0.5)
  \roundnode{n4}(11,0.5)
  \roundnode{n5}(12,0.5)
  \freetext(9.3,1){$C_3$}
  \edge{top1}{n1}
  \edge{bottom1}{n1}
  \edge{n1}{n2}
  \edge{n2}{n3}
  \autonodetext{n3}[e]{\large $\cdots$}
  \autonodetext{n4}[w]{\large $\cdots$}
  \edge{n4}{n5}
  \autonodetext{top1}[w]{$(-2)$}
  \autonodetext{bottom1}[w]{$(2)$}
  \roundnode{top3}(1,4)
  \roundnode{bottom3}(1,3)
  \roundnode{l1}(2,3.5)
  \roundnode{l2}(3,3.5)
  \roundnode{l3}(4,3.5)
  \roundnode{l4}(6,3.5)
  \roundnode{l5}(7,3.5)
  \freetext(4.3,4){$C_1$}
  \edge{top3}{l1}
  \edge{bottom3}{l1}
  \autonodetext{bottom3}[w]{$(2)$}
  \autonodetext{top3}[w]{$(-4)$}
  \edge{l1}{l2}
  \edge{l2}{l3}
  \autonodetext{l3}[e]{\large $\cdots$}
  \autonodetext{l4}[w]{\large $\cdots$}
  \edge{l4}{l5}
  \roundnode{top4}(10,4)
  \roundnode{bottom4}(10,3)
  \roundnode{p1}(11,3.5)
  \roundnode{p2}(12,3.5)
  \roundnode{p3}(13,3.5)
  \roundnode{p4}(15,3.5)
  \roundnode{p5}(16,3.5)
  \freetext(13.3,4){$C_2$}
  \edge{top4}{p1}
  \edge{bottom4}{p1}
  \autonodetext{top4}[w]{$(-4)$}
  \autonodetext{bottom4}[w]{$(2)$}
  \autonodetext{p5}[n]{$(2)$}
  \edge{p1}{p2}
  \edge{p2}{p3}
  \autonodetext{p3}[e]{\large $\cdots$}
  \autonodetext{p4}[w]{\large $\cdots$}
  \edge{p4}{p5}
\end{graph}
\end{center}
\caption{The characteristic classes determining $\t_{\xi_1}$,
$\t_{\xi_2}$ and $\t_\Xi$}
\label{f:classes} 
\end{figure}
Since when $p$ is even $2\mu_b=0$, while when $p$ is odd $\mu_c=3\mu_b$,
we have
\[
\t_{\xi_1}=\t_1-2\mu_b=
\begin{cases}
\t_1\quad\text{for $p$ even},\\
\t_1-\mu_c+\mu_b=\t_2\quad\text{for $p$ odd}
\end{cases}
\]
Since $\mu_a=\mu_c-\mu_b=\mu_b-\mu_c$, we have
\[
\t_{\xi_2}=\t_{\xi_1}+\mu_a=\t_{\xi_1}+\mu_b-\mu_c=
\begin{cases}
\t_1+\mu_b-\mu_c=\t_2\quad\text{for $p$ even},\\
\t_1\quad\text{for $p$ odd}
\end{cases}
\]
and
\[
\t_\Xi=\t_1-\mu_b=\t_1+\mu_a-\mu_c=\t_3.
\]
\end{proof}

We have shown that for $p\geq 2$ the contact structures $\xi_1$, $\xi_2$
and $\Xi$ are distinct up to homotopy. Next, we are going to prove
that they are tight for every $p\geq 2$.

\begin{defn}
Let $\eta$ be the contact structure defined by the diagram obtained
from any of the diagrams of Figure~\ref{f:structure} by omitting the
Legendrian knot $L$. 
\end{defn}

A simple Kirby calculus computation shows that $\eta$ is a contact
structure on $S^1\times S^2$.

\begin{prop}\label{p:eta}
The contact Ozsv\'ath--Szab\'o invariant of $\eta$ is nonzero.
\end{prop}

\begin{proof}
Consider the contact structure $\zeta$ given by the surgery diagram
obtained from Figure~\ref{f:structure}(a) by erasing both $L$ and one
of the $(+1)$--framed Legendrian unknots. According to
Corollary~\ref{c:legsurg} and Lemma~\ref{l:etak}, the contact
Ozsv\'ath--Szab\'o invariant of the resulting structure $\zeta$ is
nontrivial. It is easy to see that the 3--manifold underlying $\zeta$
is the lens space $L(4,1)$.  Let $X$ denote the cobordism from
$L(4,1)$ to $S^1\times S^2$ obtained by the handle attachment defined
by the remaining contact $(+1)$--surgery.  According to
Theorem~\ref{t:item2} we have
\[ 
F_{-X}(c(L(4,1), \zeta ))=c(S^1\times S^2, \eta ),
\]
where $F_{-X}=\sum_{\s\in Spin^c(X)} \pm F_{-X,\s}$. The cobordism $-X$
induces an exact triangle
\[
\begin{graph}(6,2)
\graphlinecolour{1}\grapharrowtype{2}
\textnode {A}(1,1.5){$\hf (-L(4,1))$}
\textnode {B}(5, 1.5){$\hf (S^1\times S^2)$}
\textnode {C}(3, 0){$\hf ({\mathbb {RP}}^3\# {\mathbb {RP}}^3)$}
\diredge {A}{B}[\graphlinecolour{0}]
\diredge {B}{C}[\graphlinecolour{0}]
\diredge {C}{A}[\graphlinecolour{0}]
\freetext (3,1.8){$F_{-X}$}
\freetext (4.6,0.6){$F_U$}
\end{graph}
\]
The Ozsv\'ath--Szab\'o homology groups in this triangle are well--known 
(see \cite{OSzabs}):
\begin{itemize}
\item 
$\hf (-L(4,1))\cong\Z^2_{(0)}\oplus\Z_{(\frac34)}\oplus\Z_{(-\frac14)}$
\item 
$\hf (S^1\times S^2)\cong\Z_{(\frac12)}\oplus\Z_{(-\frac12)}$ 
\item 
$\hf ({\mathbb {RP}}^3\# {\mathbb {RP}}^3)\cong
\Z^2_{(0)}\oplus\Z_{(\frac12)}\oplus\Z_{(-\frac12)}$ 
\end{itemize}
A simple computation shows that  the cobordisms $-X$ and $U$
have zero signature. Moreover, since the $4$--manifolds $X$ and $U$ are
obtained by attaching a $2$--handle to $S^1\x S^2$, the restriction
maps $H^2(X;\Z)\to H^2(\del X;\Z)$ and $H^2(U;\Z)\to H^2(\del U;\Z)$
are injective. Since the group $\hf(S^1\x S^2)$ is concentrated at the
torsion spin$^c$ structure, this implies that each possibly
nontrivial component of the maps $F_{-X}$ and $F_{U}$ is induced by
a torsion spin$^c$ structure. Therefore, by the degree--shift
formula~\eqref{e:d-shift}, both $F_{-X}$ and $F_U$ shift degrees by
$-\frac{1}{2}$.

Exactness of the triangle immediately shows that the kernel of
$F_{-X}$ is 3--dimensional and, since $F_{-X}$ shifts degree by
$-\frac12$, we see that
\[
\Z_{(\frac34)}\oplus\Z_{(-\frac14)}\subseteq \ker F_{-X}
\quad\text{and}\quad
F_{-X}(\Z^2_{(0)})=\Z_{(-\frac12)}.
\]  
Since the ${\mathcal {J}}$--map preserves degree and fixes spin
structures, the summands $\Z_{(\frac34)}$ and $\Z_{(-\frac14)}$ 
inside $\hf(-L(4,1))$ correspond to the two spin structures on $-L(4,1)$.
Therefore, since $c(\ze)\not=0$, 
\[
\langle c(\ze), \J_{-L(4,1)} c(\ze)\rangle = \Z^2_{(0)} \subseteq \hf(-L(4,1)).
\]
The group $\hf (S^1\times S^2)$ is pointwise fixed (up to sign) by the
${\mathcal J}$--action because it is concentrated at the only spin$^c$
structure induced by a spin structure, and has rank at most one in each 
degree. Hence,
\[
F_{-X}(\J_{-L(4,1)} c(\ze))=\J_{S^1\x S^2} F_{-X}(c(\ze))=
F_{-X}(c(\ze)).
\]
Since $F_{-X}(\Z^2_{(0)})\not=\{0\}$, this implies 
$ c(\eta)=F_{-X}(c(\ze))\not=0.$
\end{proof} 

\begin{cor}\label{c:stein}
The contact structures $\xi_1$ and $\xi_2$ are Stein fillable, hence
tight.
\end{cor}

\begin{proof}
The contact structure $\eta$ has nonzero contact Ozsv\'ath--Szab\'o
invariant, and therefore it is tight. It is well--known that
$S^1\times S^2$ carries a unique tight contact structure up to
isotopy, and this contact structure is Stein fillable. Therefore,
$\eta$ is Stein fillable. Since contact $(-1)$--surgery preserves
Stein fillability, the statement follows.
\end{proof}

\begin{thm}\label{t:fotight}
The Ozsv\'ath--Szab\'o invariant of the contact structure $\Xi$ is
nonzero.
\end{thm}

\begin{proof}
Denote by $(Y,\be)$ the contact 3--manifold whose contact surgery
presentation is obtained from Figure~\ref{f:bigksi} by erasing one of
the unknots with contact surgery coefficient equal to $+1$.  Let $-X$
be the cobordism from $Y$ to $M_p$ determined by the missing contact
$(+1)$--surgery, with orientation reversed. Denote by $c^+(\be)$ the
image of $c(\be)$ under the map induced by the natural
homomorphism~\cite{OSzF1, OSzF2}
\[
\hf(-Y,\t_{\be})\to 
HF^+(-Y,\t_{\be}),
\]
and define $c^+(\Xi)$ in the analogous way. Clearly, it is enough to
show that $c^+(\Xi)\not=0$.  By~\cite[Lemma~2.11]{Gh}, there is a
spin$^c$ structure $\s$ on $-X$ such that
$F^+_{-X,\s}(c^+(\be))=c^+(\Xi)$ and
\[
-d_3(\be) + \de(\s) = -d_3(\Xi),
\]
where
\[
\de(\s):=\frac 14 (c_1^2(\s) - 3\si(-X) - 2\chi(-X)).
\]

The 3--manifold $-Y$ is given by the surgery presentation obtained
from the third diagram of Figure~\ref{f:kirby} by changing the framing
of the $(1)$--framed unknot to $0$. Then, Kirby moves similar to those
of Figure~\ref{f:kirby} show that $-Y$ is the boundary of a plumbing
whose graph is obtained from the graph of Figure~\ref{f:D} by changing
the framing of the central vertex to
$-3$. By~\cite[Theorem~7.1]{OSzsymp} it follows that $-Y$ is an
$L$--space. By Lemma~\ref{l:etak} and Corollary~\ref{c:legsurg} we
have $c(\be)\not=0$, therefore $-d_3(\be) = d(-Y,\t_{\be})$. This
immediately implies $c^+(\be)\not=0$. Moreover,
using~\cite[Corollary~3.6]{DGS} (where $b_2(X)$ should be plugged into
the formula instead of the Euler characteristic $\chi(X)$ because the
3--dimensional invariant used in Heegaard Floer theory is shifted by
$\frac12$) a simple calculation gives
\[
d_3(\Xi)=\frac{2-p}4. 
\]
Therefore, by Proposition~\ref{p:mp} and Lemma~\ref{l:spin-c} we have
$-d_3(\Xi) = d(-M_p,\t_{\Xi})$.

Another simple calculation shows that $b_2^{-}(-X)=1$. Since $\t_\be$
and $\t_\Xi$ are torsion and $HF^\infty(-Y,\t_\eta)\cong\Z[U,U^{-1}]$,
by~\cite[Proposition 9.4]{OSzabs} the induced map
\[
F^{\infty}_{-X,\s}\co HF^{\infty}(-Y,\t_{\be})
\to HF^{\infty}(-M_p,\t_{\Xi})
\]
is an isomorphism. Since the group $HF^{-}_d(-Y,\t_\be)$ vanishes if
the absolute degree $d$ is sufficiently large, the commutative diagram
(see~\cite[Section~2]{OSzabs})
\[
\begin{graph}(15,2)
\graphlinecolour{1}\grapharrowtype{2}
\textnode {A1}(1,1.5){$\cdots$}
\textnode {A2}(3.5,1.5){$HF^{-}(-Y,\t_\be)$}
\textnode {A3}(7.2,1.5){$HF^{\infty}(-Y,\t_\be)$}
\textnode {A4}(11,1.5){$HF^{+}(-Y,\t_\be)$}
\textnode {A5}(13.5,1.5){$\cdots$}
\textnode {B1}(1,0){$\cdots$}
\textnode {B2}(3.5,0){$HF^{-}(-M_p,\t_\Xi)$}
\textnode {B3}(7.2,0){$HF^{\infty}(-M_p,\t_\Xi)$}
\textnode {B4}(11,0){$HF^{+}(-M_p,\t_\Xi)$}
\textnode {B5}(13.5,0){$\cdots$}
\diredge {A1}{A2}[\graphlinecolour{0}]
\diredge {A2}{A3}[\graphlinecolour{0}]
\diredge {A3}{A4}[\graphlinecolour{0}]
\diredge {A4}{A5}[\graphlinecolour{0}]
\diredge {B1}{B2}[\graphlinecolour{0}]
\diredge {B2}{B3}[\graphlinecolour{0}]
\diredge {B3}{B4}[\graphlinecolour{0}]
\diredge {B4}{B5}[\graphlinecolour{0}]
\diredge {A4}{B4}[\graphlinecolour{0}]
\diredge {A2}{B2}[\graphlinecolour{0}]
\diredge {A3}{B3}[\graphlinecolour{0}]
\freetext (4.3,0.8){$F^-_{-X,\s}$}
\freetext (8,0.8){$F^\infty_{-X,\s}$}
\freetext (11.7,0.8){$F^+_{-X,\s}$}
\end{graph}
\]
implies that the map 
\[
F^+_{-X,\s}\co HF^+_d (-Y,\t_{\be})
\to HF^+_{d+\de(\s)} (-M_p,\t_{\Xi})
\]
is also an isomorphism when $d$ is large enough. Since $F^+_{-X,\s}$
is a homomorphism of $\Z[U]$--modules~\cite{OSzF4}, this immediately
implies that $F^+_{-X,\s}$ restricted to
$HF^+_{-d_3(\be)}(-Y,\t_{\be})$ is an isomorphism if and only if
$-d_3(\Xi) = d(-M_p,\t_{\Xi})$, and the conclusion follows.
\end{proof}

\begin{rem}\label{r:more}
(1) The proof of Corollary~\ref{c:stein} applies to show that the result
of any Legendrian surgery on $(S^1\x S^2,\eta)$ is Stein fillable. In
particular, for any choice of zig--zag distribution for the Legendrian
knot $L$ of Figure~\ref{f:structure}, the resulting contact structures
$\xi _i $ ($i=1,2, \ldots , p$) are Stein fillable. 
In addition, by reversing the stabilizations on the other two
contact $(-1)$--framed knots of Figure~\ref{f:structure}, another
collection of tight contact surgery diagrams --- denoted by $\xi _i '$ 
($i=1,2,\ldots, p $) --- can be given. Although these diagrams give isotopic 
structures to $\xi _1$ or $\xi _2$ on $M_p$, they will play 
an important role in the classification results discussed in the next 
section.

(2) Let $\Xi'$ be the contact structure on $M_p$ with surgery
presentation obtained from Figure~\ref{f:bigksi} by applying a
$180^{\circ}$ rotation around an axis perpendicular to the plane of
the picture. Then, $\t_{\Xi'}=\t_{\Xi}$ and
$d_3(\Xi')=d_3(\Xi)$. Since the auxiliary 3--manifold $Y$ used in the
proof of Theorem~\ref{t:fotight} is the same for $\Xi'$, the same
proof also shows that $\Xi'$ has nonzero Ozsv\'ath--Szab\'o invariant.

(3) For a more general form of Theorem~\ref{t:fotight} see
    \cite{OSzIII}.
\end{rem}

The results above lead to
 
\begin{cor}\label{c:mpclass}
For every $p\geq 2$, the 3--manifold $M_p=M(-1;
\frac{1}{2},\frac{1}{2},\frac{1}{p})$ admits exactly three tight
contact structures (up to isotopy).
\end{cor}

\begin{proof}
Corollaries~\ref{c:ub1} and~\ref{c:ub2} imply that $M_p$ admits at
most three tight contact structures, while to combination of
Lemma~\ref{l:spin-c}, Corollary~\ref{c:stein} and
Theorem~\ref{t:fotight} verifies that $M_p$ admits at least three
distinct tight contact structures, concluding the proof. Notice that
this argument shows, for example, that $\Xi$ and $\Xi'$ are isotopic
on $M_p$.
\end{proof}

\sh{Nonfillability of $(M_p,\Xi)$}

In this subsection we give simple proofs of the facts that the tight
contact 3--manifold $(M_p,\Xi)$ is not Stein fillable for $p>2$, and
not symplectically fillable for $p\not\equiv 2\bmod 8$, justifying our use of
contact Ozsv\'ath--Szab\'o invariants in the proof of their
tightness. First we need the following general observation:

\begin{lem} \label{l:fill}
Let $(Y,\xi)$ be a contact 3--manifold such that $b_1(Y)=0$, and
suppose that for every symplectic filling $(X,\om)$ of $(Y,\xi)$ we
have $b_2^+(X)=0$. Then, for every symplectic filling $(X,\om)$ of
$(Y,\xi)$ we also have $b_1(X)=0$.
\end{lem}

\begin{proof}
Let $n$ be the order of the finite group $H_1(Y;\Z)$. By
contradiction, suppose that $(X,\om)$ is a symplectic filling of
$(Y,\xi)$ with $b_1(X)>0$. Then $X$ admits a connected $(n+1)$--fold
cover $\tilde X$ which is necessarily trivial over $\partial X
=Y$. By~\cite{El, Et} we can cap off $n$ of the boundary components of
$\tilde X$ with symplectic caps having $b_2^+>0$, obtaining a
symplectic filling of $(Y,\xi)$ with $b_2^+>0$, which is against our
assumptions.
\end{proof}

\begin{thm}\label{t:nonfill}
The tight contact 3--manifold $(M_p, \Xi )$ is not 
Stein fillable for every $p>2$, and is not symplectically fillable
for $p\not\equiv 2\bmod 8$.
\end{thm}

\begin{proof}
By~\cite[Theorem~2.2]{Ligok} if $(X, \omega)$ is any symplectic
filling of $(M_p,\Xi)$, then the intersection form $Q_X$ is a standard
diagonal negative definite form. Therefore, Lemma~\ref{l:fill} implies
$b_1(X)=0$.

Now let $(X,J)$ be a Stein filling of $(M_p,\Xi)$ with complex
structure $J$ inducing the contact structure $\Xi$. Then, the complex
structure $-J$ gives another Stein structure on $X$ inducing a contact
structure $\overline\Xi$ on $M_p$. Since the associated spin$^c$
structure satisfies $\s_{-J}=\overline{\s_J}$, we have 
\[
\t_{\overline\Xi}=\overline{\t_\Xi}=\overline{\t_3}=\t_3
\]
because, as we observed previously, $\t_3$ is induced by a spin
structure. Therefore, by Corollary~\ref{c:mpclass} $\overline{\Xi}$ is
isotopic to $\Xi$.  But then, by~\cite[Theorem~1.2]{LM} we have
$\overline{\s_J}=\s_J$, which implies $c_1(J)=0$. Since $Q_X$ is
standard diagonal and $b_1(X)=0$, this implies $\si(X)=0$ and
$\chi(X)=1$. In view of Formula~\eqref{e:d3-inv} we have $d_3(\Xi)=0$, and by 
Proposition~\ref{p:mp} this is possible only if
$p=2$.

When $(X,\om)$ is a general symplectic filling of $(M_p,\Xi)$ we are
unable to conclude that $c_1(X, \om)=0$, but we still know that $Q_X$ is
negative definite, diagonal and $b_1(X)=0$, therefore 
Formula~\eqref{e:d3-inv} and Proposition~\ref{p:mp} give
\[
\frac{2-p}{4}=\frac{1}{4}(-\sum _{i=1}^{b_2(X)} (2n_i+1)^2+b_2(X))
\]
for some $n_i\geq 0$. Since $(2n_i+1)^2-1=4n_i(n_i+1)$ is divisible by
8, the equation provides the desired contradiction once $(p-2)$ is not
divisible by $8$.
\end{proof}

It is natural to expect that $(M_p,\Xi)$ is not symplectically
fillable for every $p>2$. On the other hand, from the obvious $\Z
/3\Z$--symmetry of the surgery diagram of $M_2=M(-1; \frac{1}{2},
\frac{1}{2}, \frac{1}{2})$ it is not hard to see that for $p=2$ the
structures $\xi _1, \xi _2$ and $\Xi$ are contactomorphic (although
not isotopic), hence for $p=2$ the structure $\Xi$ is Stein fillable.

\section{Lower bounds and the proof of Theorem~\ref{t:main}}
\label{s:general}

Now we are ready to prove the general lower bounds for the number of
tight contact structures on the manifolds under consideration.  We
will proceed by first constructing a set of contact structures which
are --- due to our previous computations --- all tight, and then
determining how many distinct structures are in that set. Consider the
surgery presentations of the contact structures $\xi _i$, $\xi _i '$
($i=1,2,\ldots , p$), $\Xi$ and $\Xi '$ on $M_p$
(cf. Figure~\ref{f:structure}, Figure~\ref{f:bigksi} and
Remarks~\ref{r:more}(1) and (2)). According to the continued fraction
expansions of the surgery coefficients $-\frac{1}{r_i}$, attach chains
of Legendrian unknots $K_j ^i$ $(i=1,2,3; j=1,\ldots , k_i)$
stabilized $(a_j^i -2)$--times to the contact $(-1)$--framed knots of
Figures~\ref{f:structure} and \ref{f:bigksi}. Notice that there are
many choices for the required stabilizations, hence this procedure
gives rise to a number of contact structures.

Define $A(\xi )$ as the set of contact structures on $M(-1; r_1, r_2,
r_3)$ obtained by Legendrian surgery on the knots $K_j^i$ on either
diagram of Figure~\ref{f:structure} or its modifications $\xi _i , \xi
_i '$ ($i=1,2, \ldots , p$) described in Remark~\ref{r:more}(1).  In a
similar manner, $A(\Xi )$ denotes the set of contact structures
obtained either from the diagram of Figure~\ref{f:bigksi} giving $\Xi$
or its symmetric giving $\Xi'$.

\begin{lem}
The set $A(\xi)\cup A(\Xi)$ consists of tight contact structures
having nonzero Ozsv\'ath--Szab\'o invariant.
\end{lem}

\begin{proof}
Any element of $A(\xi)\cup A(\Xi)$ is constructed by Legendrian
surgery on a contact structure with nonzero contact Ozsv\'ath--Szab\'o
invariant, therefore the statement follows immediately from
Corollary~\ref{c:legsurg}.
\end{proof}

\begin{prop}\label{p:disjoint}
If $\zeta_1\in A(\xi)$ and $\zeta_2 \in A(\Xi)$ then
$c(M,\zeta_1)\neq c(M,\zeta_2)$. In particular, $\ze_1$ is not
isotopic to $\ze_2$.
\end{prop}

\begin{proof}
Denote by $(Y_p,\ga_i)$, $i=1,2$ the contact 3--manifolds obtained by
contact $(+1)$--surgeries in $(M,\ze_i)$ along a push--off of
$K_1^i$ for every $i$. The 3--manifold $Y_p$ is a connected sum
$M_p\# L$, where $L$ is the connected sum of at most three lens
spaces. Correspondingly, the contact structure $\ga_i$ can be written
as
\[
\ga_i=\ga_i^{M_p}\# \ga_i^L,
\]
where $\ga_1^{M_p}$ is equal to $\xi_i$ or $\xi_i '$ ($i=1,2, \ldots ,
p$) and $\ga_2^{M_p}$ to $\Xi$ or $\Xi'$. By Corollary~\ref{c:stein}
and Theorem~\ref{t:fotight} $c(M_p, \ga_1^{M_p})$ and $c(M_p,
\ga_2^{M_p})$ are both nonzero. Since by Lemma~\ref{l:spin-c} $c(M_p,
\ga_1^{M_p})$ and $c(M_p, \ga_2^{M_p})$ live in groups corresponding
to different spin$^c$ structures, we have $c(\ga_1^M)\neq
c(\ga_2^M)$. Thus, the map corresponding to the cobordism induced by
the prescribed contact $(+1)$--surgeries sends $c(M, \zeta_1)$ and
$c(M,\zeta_2)$ to distinct elements, and the statement follows.
\end{proof}

As a consequence of Proposition~\ref{p:disjoint}, in order to get a lower
bound for the number of nonisotopic contact structures on $M$, we can
examine the sets $A(\xi )$ and $A(\Xi )$ separately.

\sh{Lower bound on $|A(\xi)|$}

\begin{figure}[ht]
\begin{center}
\psfrag{-l+1}{$-l+1$}
\psfrag{-p+2}{$-p+2$}
\psfrag{-k}{$-k$}
\psfrag{2}{$2$}
\psfrag{(a)}{$(a)$}
\psfrag{(b)}{$(b)$}
\psfrag{-l}{$-l$}
\psfrag{-2}{$-2$}
\psfrag{-p}{$-p$}
\psfrag{-1}{$-1$}
\includegraphics[width=11cm]{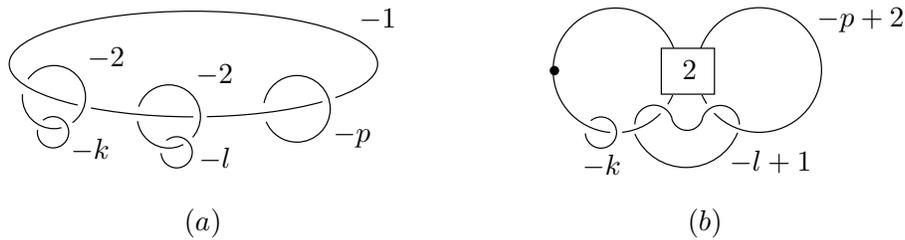}
\end{center}
\caption{\quad Surgery diagrams for (a) $M_{p,k,l}$ and (b) $X_{p,k,l}$}
\label{f:xpkl}
\end{figure}

Suppose first that $k_1=k_2=1$ and $k_3=0$, that is, there are
two circles on the first and second legs and there is a single
one on the third. Let the corresponding 3--manifold be denoted by
$M_{p,k,l}$. A surgery presentation for this 3--manifold is given by
Figure~\ref{f:xpkl}(a). Notice that if $r_1>r_2=\frac12$ then
$k_2=0$. Therefore, to cover this case we shall also consider the
3--manifold $M_{p,k}$ defined as in Figure~\ref{f:xpkl}(a) but
omitting the $(-l)$--framed knot. (The case $r_1=r_2=\frac12$ leads to
the manifold $M_p$--- we have already dealt with this manifold in
Section~\ref{s:special}.)

\begin{prop}\label{p:count}
On the 3--manifold $M_{p,k,l}$ there are at least
\[
(2(k-1)(l-1) + p-1)(k+l-2)
\]
isotopy classes of tight contact structures belonging to $A(\xi)$, while 
on $M_{p,k}$ there are at least 
\[
2(k-1) + p-1
\]
such isotopy classes.
\end{prop}

\begin{proof}
Let us start with the case of $M_{p,k,l}$, i.e.~when $r_2>\frac12$. 
We will apply~\cite[Theorem~1.2]{LM}, which implies that if two Stein
structures on a 4--manifold $X$ have distinct first Chern classes,
then the induced contact structures on $\partial X$ are
nonisotopic. As proved by Plamenevskaya~\cite{Olga}, tight contact
structures distinguished in this way have different contact
Ozsv\'ath--Szab\'o invariants. Notice, however, that the contact
surgery diagrams giving the elements of $A(\xi)$ do not provide Stein
fillings. A simple surgery operation, however, can turn the
4--manifold $W_{p,k,l}$ given by each surgery diagram into a Stein
domain.  Namely, let us consider the codimension--0 submanifold
$Z\subset W_{p,k,l}$ defined by the union of the two $(+1)$--framed
Legendrian unknots together with the two once stabilized
$(-1)$--knots. By Proposition~\ref{p:eta}, the corresponding contact
structure $\eta$ is the unique tight (and hence Stein fillable)
contact structure on $S^1\x S^2$.

Replacing $Z$ with a 4--dimensional 1--handle $H$ we obtain a
4--manifold
\[
X_{p,k,l}=(W_{p,k,l}\setminus Z)\cup H
\]
with a decomposition involving a 1--handle and three 2--handles. Then,
$X_{p,k,l}$ can be thought of as obtained by attaching three Stein
2--handles to a Stein 1--handle, and therefore carries a Stein
structure. The Legendrian attaching circles are $L$ from
Figure~\ref{f:structure} plus two Legendrian meridional unknots $M_1$
and $M_2$ linking the once stabilized unknots in the same picture.
Smoothly, a handlebody decomposition for $X_{p,k,l}$ is given by
Figure~\ref{f:xpkl}(b), where the $(2-p)$--framed knot corresponds to
$L$, and $M_1$, $M_2$ correspond, respectively, to the $(-k)$-- and
the $(1-l)$--framed knots. Suppose now that (for some choice of the
orientations) the rotation numbers of the once stabilized Legendrian
unknots in Figure~\ref{f:structure} are $A$ and $-A$ with $A\in\{\pm
1\}$, while the rotation numbers of $M_1$, $M_2$ and $L$ are given,
respectively, by $x$, $y$ and $z$.  These rotation numbers
satisfy the constraints
\begin{equation}\label{e:constraints}
\begin{split}
x & \in\{-k+2i_1,\ i_1=1,\ldots, k-1\},\\ 
y & \in\{-l+2i_2,\ i_2=1,\ldots,l-1\},\\
z & \in\{-p-1+2i_3,\ i_3=1,\ldots,p\}.
\end{split}
\end{equation}
(Note the special behaviour of $z$, which is the rotation number of
$L$, linking the $(+1)$--surgery curves in Figure~\ref{f:structure}.)
Denote by $a$, $b$ and $c$ the homology classes in $H_2(W_{p,k,l};\Z)$
defined by $M_1$, $M_2$ and $L$ respectively, and observe that there
are homology classes $\al,\be\in H_2(W_{p,k,l}\setminus Z;\Z)$ which
map to $c-a-b$ and $a-b$, respectively, and such that their images
$\overline\al$ and $\overline\be$ under the map induced by the
inclusion $(W_{p,k,l}\setminus Z)\subset X_{p,k,l}$ generate the group
$H_2(X_{p,k,l};\Z)\cong\Z^2$. It is not hard to see that $X_{p,k,l}$
is simply connected, and the values of the first Chern class of its
Stein structure on $\overline\al$ and $\overline\be$ are,
respectively, $z-x-y+2$ and $x-y-A+2$. Therefore, to
apply~\cite[Theorem~1.2]{LM} we need to count the number of elements
of the set $S(p,k,l)\subset\Z^2$ of pairs $(x-y-A,z-x-y)$ such that
$A\in\{\pm 1\}$ and $x,y,z$ satisfy the constraints of
\eqref{e:constraints}.

In order to do this, we first consider the set $T(k,l)$ consisting of
pairs $(x-y,-x-y)$, where $x$ and $y$ satisfy the constraints given
by~\eqref{e:constraints}. Setting $e_1=\binom10$ and $e_2=\binom01$, we 
have 
\[
S(p,k,l)=\pm e_1 + \bigcup_z \left(T(k,l)+ze_2 \right) ,
\]
where $z$ satifies~\eqref{e:constraints}. Let $\varphi\co\Z^2\to\Z^2$
be the injective map given by $\varphi(x,y)=(x-y,-x-y)$. Clearly
$T(k,l)=\varphi(Q(k,l))$, where $Q(k,l)\subset\Z^2$ is the set of
pairs $(x,y)$ which satisfy the constraints given
by~\eqref{e:constraints}. Clearly, $Q(k,l)$ has the shape of a
rectangle and contains $(k-1)(l-1)$ elements. 
Since 
\[
\varphi(\tfrac12,-\tfrac12)=e_1,\quad
\varphi(-\tfrac12,-\tfrac12)=e_2,
\]
the set $S(p,k,l)$ has the same number of elements as the set 
\[
T(p,k,l):=
\pm\begin{pmatrix}\phantom{-}1/2\\-1/2\end{pmatrix} + 
\bigcup_z\left(Q(k,l)+z\begin{pmatrix}-1/2\\-1/2\end{pmatrix}\right)
\subset\Z^2. 
\]
To count the number of elements in this set, observe that the
$(k-1)(l-1)$ elements of $Q(k,l)$ form a rectangle in the plane, and
they are at two integral units of distance from each other.  In
Figure~\ref{f:Q} the points of the set $Q(5,4)$ are represented by
a's.
\begin{figure}[ht]
\begin{center}
\psfrag{a}{a}
\psfrag{b}{b}
\psfrag{c}{c}
\includegraphics[height=4.5cm]{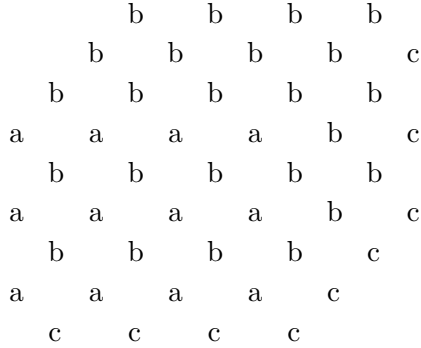}
\end{center}
\caption{Counting the number of elements of $T(p,k,l)$}
\label{f:Q}
\end{figure}
The set 
\[
R(p,k,l):=\bigcup_z\left(Q(k,l)+
z\begin{pmatrix}-1/2\\-1/2\end{pmatrix}\right)
\]
is obtained as the union of $p$ shifts of $Q(k,l)$ by integral units
in the North--East direction.  If we assume $p=4$, for instance, we
see in Figure~\ref{f:Q} how the results of these shifts create the
elements denoted by b's in the picture. It is easy to compute that the
number of elements increases by $(k-1)(l-1)+(k+l-2)(p-2)$. Finally,
$T(p,k,l)$ is obtained as the union of $2$ shifts of $R(p,k,l)$, one 
integral unit apart from each other in the South--East
direction. In Figure~\ref{f:Q} the resulting new elements have been
denoted by c's. It is easy to see that the number of such elements is
$(k-1) + (l-1) + (p-2)$. Therefore, the cardinality of $T(p,k,l)$ is
obtained by adding the number of a's, b's and c's:
\[
\begin{split}
(k-1)(l-1) + (k-1)(l-1)+(k+l-3)(p-2) + (k-1) + (l-1) + (p-2) = \\
2(k-1)(l-1) + (p-1)(k+l-2).
\end{split}
\]
When $r_1>r_2=\frac12$, i.e.~in the case of $M_{p,k}$, there is no
meridian $M_2$, no variable $y$ nor homology class $b$, and one can
work with 4--manifolds $W_{p,k}$ and $X_{p,k}$ by analogy to what we
did before. There is a class $\al\in H_2(W_{p,k}\setminus Z;\Z)$ which
is sent to $c-2a\in H_2(W_{p,k};\Z)$ by the map induced by inclusion,
and whose image $\overline\al\in H_2(X_{p,k};\Z)\cong\Z$ is a
generator. $X_{p,k}$ is still simply connected, and the possible
values of the first Chern classes of $X_{p,k}$ on $\overline\al$ are
of the form $z+A-2x$ while $z$, $A$ and $x$ range as
in~\eqref{e:constraints}. An easy count yields the stated formula.
\end{proof}

\begin{cor}\label{c:Axi}
The number of isotopy classes of tight contact structures on the 
3--manifold $M(-1;r_1,r_2,r_3)$ belonging to $A(\xi)$ are at least
\[
\left(2(a_1^1-1)(a_1^2-1) + (a_0^3-1)(a_1^1+a_1^2-2)\right)
(a^3_1-1)
\prod_{i=1}^3 \prod_{j\geq 2}(a_j^i-1)
\]
if $r_2>\frac12$,
\[
\left(2(a_1^1 -1) + (a_0^3 -1)\right)(a_1^3-1)\prod_{i\neq
2}\prod_{j\geq 2}(a^i_j-1)
\]
if $r_1>r_2=\frac 12$, and 
\[
2\prod_{j\geq 1}(a_j^3-1)
\]
if $r_1=r_2=\frac12$. In the above formulae $a_j^i=2$ by convention
if $j>k_i$.
\end{cor}

\begin{proof}
Let $\ze$ be an element of the set $A(\xi)$. Perform contact
$(+1)$--surgeries along the Legendrian push--offs of $K_1^1$ (if it
exists), $K_1^2$ (if it exists) and $K_0^3$. By
Proposition~\ref{p:eta}, the resulting contact 3--manifold is the
tight contact $S^1\times S^2$ connected sum with at most three tight
contact lens spaces. But for such contact lens spaces it is known that
the zig--zag distribution is determined by the contact invariant.
Therefore, applying Theorem~\ref{t:item2} we see that different
zig--zag distributions in the diagram for $\ze$ after the second
circle on the first two legs and after the first circle on the third
leg yield nonisotopic structures. If $r_2>\frac12$ or $r_1>r_2=\frac12$ 
the statement follows by Proposition~\ref{p:count} and a simple
computation. If $r_1=r_2=\frac12$, by Lemma~\ref{l:spin-c}
and Corollary~\ref{c:stein}, the set $A(\xi)$ contains at least $2$
elements in the case of the 3--manifold $M_p$. Thus, the stated
formula follows immediately.
\end{proof}

\sh{Lower bound on $|A(\Xi)|$ and the proof of Theorem~\ref{t:main}}  

\begin{figure}[ht]
\begin{center}
\psfrag{-2}{$-2$}
\psfrag{-k}{$-k$}
\psfrag{-l}{$-l$}
\psfrag{-p}{$-p$}
\psfrag{-1}{$-1$}
\psfrag{-m}{$-m$}
\includegraphics[height=3.5cm]{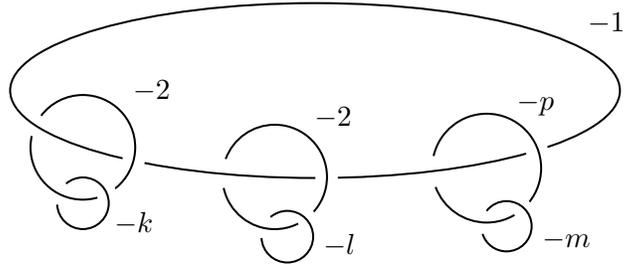}
\end{center}
\caption{\quad Kirby diagram for the 3--manifold $M_{p,k,l,m}$}
\label{f:rajz}
\end{figure}

\begin{figure}[ht]
\begin{center}
\psfrag{K1}{\small $K_1$}
\psfrag{K2}{\small $K_2$}
\psfrag{K3}{\small $K_3$}
\psfrag{K4}{\small $K_4$}
\psfrag{K5}{\small $K_5$}
\psfrag{K6}{\small $K_6$}
\psfrag{K7}{\small $K_7$}
\psfrag{K8}{\small $K_8$}
\psfrag{p-1}{\small $p-1$}
\psfrag{cusps}{\small left cusps}
\includegraphics[width=0.8\textwidth]{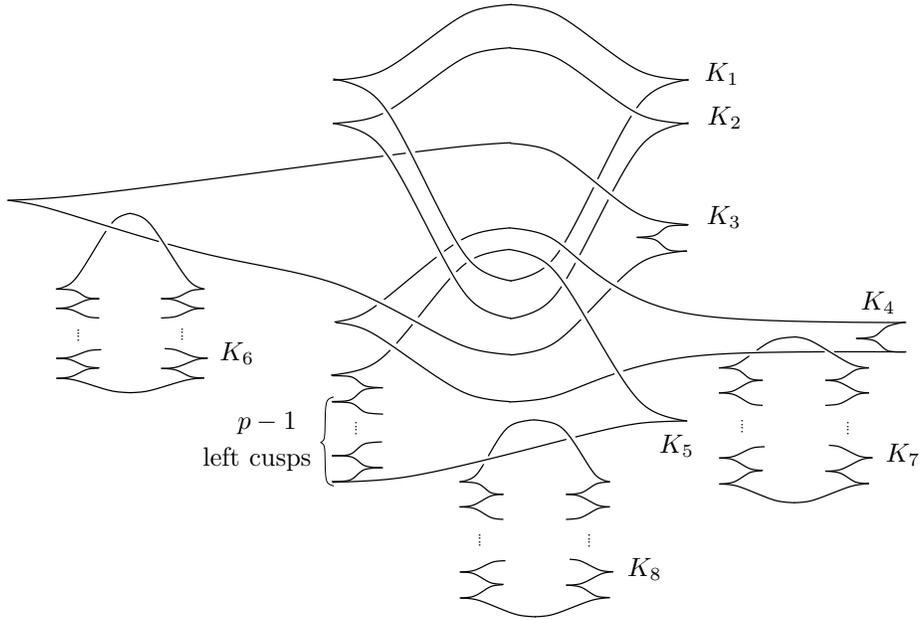}
\end{center}
\caption{\quad Contact surgery diagram for $\ze$ on $M_{p,k,l,m}$}
\label{f:Mpklm}
\end{figure}

If $k_1=k_2=k_3=1$, the same idea used to study the set
$A(\xi)$ suggests the existence of an appropriate function $g(p,k,l,m)$
such that there are at least
\[
g(a_0^3,a_1^1, a_1^2,a_1^3)  
\cdot  \Pi _{i=1}^3 \Pi _{j\geq 2}( a_j ^i -1)
\]
different elements in $A(\Xi)$: just perform contact $(+1)$--surgeries
along the push--offs of the Legendrian curves $K^i_1$ ($i=1,2,3$). (It
turns out that it is not useful to do surgery along the push--off of
the first circle of the third leg, because in the present case the
resulting contact structure on $S^1\times S^2$ would be overtwisted.)
So our aim will be to find a lower bound for the number of distinct
structures in $A(\Xi)$ on the 3--manifolds $M_{p,k,l,m}$ defined by
Figure~\ref{f:rajz}. To cover the cases when $k_i=0$ for some
$i\in\{1,2,3\}$, we shall consider also analogously defined manifolds
$M_{p,k,m}$, $M_{p,k,l}$ and $M_{p,k}$.
Let $K_1, \ldots ,K_5$ be the components of the contact surgery
diagram in Figure~\ref{f:bigksi} defining $M_p$ with one of the two
tight contact structures $\Xi$ or $\Xi'$.  Let $K_6, K_7, K_8$ be the
three extra knots linked once to $K_3, K_4, K_5$ respectively as shown
in Figure~\ref{f:Mpklm}, which gives a contact surgery presentation of
$(M_{p,k,l,m},\ze)$ with $\ze\in A(\Xi)$.
Since the contact surgery coefficient of $K_6$, $K_7$ and $K_8$ is
$-1$, they determine a Stein cobordism $W_{p,k,l,m}$ between
$(M_p,\Xi)$ or $(M_p,\Xi')$ and $(M_{p,k,l,m},\ze)$. Denote by $\t$
the spin$^c$ structure induced on $W_{p,k,l,m}$ by the Stein
structure. The contact surgery diagram of Figure~\ref{f:Mpklm}
determines also a 4--manifold $X_{p,k,l,m}$ bounded by $M_{p,k,l,m}$
and a spin$^c$ structure $\s$ on $X_{p,k,l,m}$. Let $X_p$ be the
4--manifold bounding $M_p$ obtained by surgery on the link in
Figure~\ref{f:bigksi}. Since this link is a sublink of the link in
Figure~\ref{f:Mpklm}, $X_p$ is a submanifold of $X_{p,k,l,m}$ and
$W_{p,k,l,m}= X_{p,k,l,m} \setminus X_p$ is the above mentioned
cobordism between $M_p$ and $M_{p,k,l,m}$. Moreover,
\[
\s |_{W_{p,k,l,m}}=\t.
\]
The above discussion remains essentially unchanged if the knot
$K_7$, the knot $K_8$ or both the knots $K_7$ and $K_8$ are omitted
from Figure~\ref{f:Mpklm}. In fact, it suffices to replace the triple 
\[
(M_{p,k,l,m},W_{p,k,l,m},X_{p,k,l,m}) 
\]
by, respectively, the triples 
\[
(M_{p,k,m},W_{p,k,m},X_{p,k,m}), \qquad (M_{p,k,l},W_{p,k,l},X_{p,k,l})
\qquad {\mbox {and}} \qquad (M_{p,k},W_{p,k},X_{p,k}).
\]

\begin{lem}~\label{l:goofy} 
Consider two contact surgery diagrams as in Figure~\ref{f:Mpklm},
where in both diagrams $K_7$, $K_8$ or both might be missing.  Denote
by $\zeta_1$ and $\zeta_2$ the tight contact structures induced on
$M_{p,k,l,m}$, $M_{p,k,m}$, $M_{p,k,l}$ or $M_{p,k}$, and by $\t _1$
and $\t _2$ the corresponding spin$^c$ structures induced respectively
on $W_{p,k,l,m}$, $W_{p,k,m}$, $W_{p,k,l}$ or $W_{p,k}$. If $\zeta_1$
is isotopic to $\zeta_2$ then $\t _1$ is isomorphic to $\t _2$.
\end{lem}

\begin{proof}
Consider the case of $M_{p,k,l,m}$. If $\zeta_1$ is isotopic to
$\zeta_2$ then $c(\zeta_1)=c(\zeta_2)$.  By \cite[Lemma~2.11]{Gh}
$F_{\overline{W}_{p,k,l,m}, \s}(c(\zeta_i))=c(\Xi)\neq 0$ if $\s = \t _i$
and $F_{\overline{W}_{p,k,l,m}, \s}(c(\zeta_i))=0$ for any other
spin$^c$ structure $\s$ on $W_{p,k,l,m}$, where
$\overline{W}_{p,k,l,m}$ denotes the cobordism $W_{p,k,l,m}$ viewed
upside down. The same argument applies to $W_{p,k,m}$, $W_{p,k,l}$ and
$W_{p,k}$. This immediately implies the statement.
\end{proof}

\begin{lem}\label{l:ziopaperone}
Let $\zeta_1$, $\zeta_2 \in A(\Xi)$ be two contact structures on
$M_{p,k,l,m}$, $M_{p,k,m}$, $M_{p,k,l}$ or $M_{p,k}$ given by contact
surgery diagrams as in Figure~\ref{f:Mpklm}, and let $\t _1$ and $\t
_2$ be the spin$^c$ structures on the cobordism $W_{p,k,l,m}$,
$W_{p,k,m}$, $W_{p,k,l}$ or $W_{p,k}$ induced by the contact surgery
diagrams. Denote by $x_i$, $y_i$, $z_i$, $i=1,2$, respectively, the
rotation numbers (for some choice of orientations) of the Legendrian
knots $K_6$, $K_7$ and $K_8$ for $\ze_1$ and $\ze_2$. Then,
$\t _1$ is isomorphic to $\t _2$ if and only if one of the following
conditions hold:
\begin{enumerate}
\item $\zeta_1$ and $\zeta_2$ are both built by Legendrian surgery on
$\Xi$ or both on $\Xi'$, and 
\[
(x_1,y_1,z_1)=(x_2,y_2,z_2)\quad\text{for}\quad
M_{p,k,l,m},
\]
\[
(x_1,z_1)=(x_2,z_2)\quad\text{for}\quad
M_{p,k,m},
\]
\[
(x_1,y_1)=(x_2,y_2)\quad\text{for}\quad
M_{p,k,l},
\]
\[
x_1=x_2\quad\text{for}\quad
M_{p,k}
\]
\item $\zeta_1$ is built by Legendrian surgery on $\Xi$ and $\zeta_2$
is built by Legendrian surgery on $\Xi'$, and 
\[
(x_1,y_1,z_1)=(x_2,y_2,z_2-2)\quad\text{for}\quad
M_{p,k,l,m},
\]
\[
(x_1,z_1)=(x_2,z_2-2)\quad\text{for}\quad
M_{p,k,m},
\]
\[
(x_1,y_1)=(x_2,y_2)\quad\text{for}\quad
M_{p,k,l},
\]
\[
x_1=x_2\quad\text{for}\quad
M_{p,k}
\]
\end{enumerate}
\end{lem}

\begin{proof} 
We consider first the case of $M_{p,k,l,m}$. Associated to any knot
$K_i$ in the contact surgery diagram in Figure~\ref{f:Mpklm} there is
a surface $\Sigma_i \subset X_{p,k,l,m}$ obtained by capping off a
Seifert surface of $K_i$ with the core of the $2$--handle attached
along $K_i$. The homology classes represented by the surfaces
$\Sigma_i$ freely generate $H_2(X_{p,k,l,m}; \Z)$. Denote by
$[\Sigma_1]^*, \ldots ,[\Sigma_8]^*$ the dual basis of
$H^2(X_{p,k,l,m}; \Z)$. The meridional discs $N_i$ of $K_i$ represent
relative homology classes which freely generate $H_2(X_{p,k,l,m},
M_{p,k,l,m}; \Z)$. Denote by $[N_1]^*, \ldots ,[N_8]^*$ the dual basis
of $H^2(X_{p,k,l,m}, M_{p,k,l,m}; \Z)$. The cohomology exact sequence
for the pair $(X_{p,k,l,m}, W_{p,k,l,m})$ together with the excision
isomorphism
\[
H^2(X_{p,k,l,m}, W_{p,k,l,m}; \Z) \cong H^2(X_p, M_p; \Z)
\]
gives the short exact sequence
\[
0\lra H^2(X_p, M_p; \Z)\stackrel{\varphi^*}\lra
H^2(X_{p,k,l,m}; \Z)\lra H^2(W_{p,k,l,m}; \Z)\lra 0
\]
where the map $\varphi^*$ is defined for $i=1, \ldots ,5$ as
$$\varphi^*([N_i]^*)= \sum_{j=0}^8 {\tt lk}(K_i, K_j)[\Sigma_j]^*$$
where ${\tt lk}(K_i, K_j)$ denotes the linking number between $K_i$
and $K_j$ if $i \neq j$, and the smooth surgery coefficient of $K_i$
if $i=j$. In terms of the dual bases chosen above, the map $\varphi^*$
is given by the matrix
\[
\Phi^*= \begin{pmatrix}
 0 & -1 & -1 & -1 & -1   \\
-1 &  0 & -1 & -1 & -1   \\
-1 & -1 & -3 & -1 & -1   \\
-1 & -1 & -1 & -3 & -1   \\
-1 & -1 & -1 & -1 & -p-1 \\
 0 &  0 & -1 &  0 &  0   \\
 0 &  0 &  0 & -1 &  0   \\
 0 &  0 &  0 &  0 & -1   
\end{pmatrix} 
\]
Let $\s $ be a spin$^c$ structure on $X_{p,k,l,m}$ defined by
contact surgery on the Legendrian link $K_1 \cup \ldots \cup K_8$
describing $\zeta$. By~\cite[Proposition~2.3]{G}
and~\cite[Proposition~3.4]{DGS} the first Chern class of
$\s$ is given by the formula
\[
c_1(\s )= \sum_{i=0}^8 {\tt rot}(K_i)[\Sigma_i]^*, 
\]
where ${\tt rot}(K_i)$ denotes the rotation number of the Legendrian
knot $K_i$. Since $K_1$ and $K_2$ are Legendrian unknots with
Thurston--Bennequin invariant ${\tt tb}(K_1)= {\tt tb}(K_2)= -1$,
their rotation numbers are ${\tt rot}(K_1) = {\tt rot}(K_2) =0$. If
$K_1 \cup \ldots \cup K_5$ is a contact surgery diagram for $\Xi$,
then (for a suitable choice of orientations) ${\tt rot}(K_3) = {\tt
rot}(K_4) = +1$ and ${\tt rot}(K_6)=-(p-1)$.  If it is a contact
surgery diagram for $\Xi'$ then ${\tt rot}(K_3) = {\tt
rot}(K_4) = -1$, and ${\tt rot}(K_6) =(p-1)$.

Consider two contact surgery diagrams describing tight contact
structures $\zeta_1$, $\zeta_2 \in A(\Xi)$ on $M_{p,k,l,m}$ and
inducing spin$^c$ structures $\s _1$ and $\s _2$ on
$X_{p,k,l,m}$. Since $X_{p,k,l,m}$ is simply connected, the
restrictions $\t _1$ and $\t _2$ of $\s _1$ and $\s _2$ to
$W_{p,k,l,m}$ are isomorphic if and only if
$$\frac 12 (c_1(\s _1)-c_1(\s _2))|_{W_{p,k,l,m}}=0.$$ 
If $\zeta_1$ and $\zeta_2$ are both built from $\Xi$ or
from $\Xi'$, then
\begin{equation}\label{eq:1}
\frac 12 (c_1(\s _1) - c_1(\s _2))= \frac 12 ((x_1
- x_2) [\Sigma_6]^* + (y_1 - y_2) [\Sigma_7]^* + (z_1 - z_2)
[\Sigma_8]^*),
\end{equation}
while if $\zeta_1$ is built from $\Xi$ and $\zeta_2$ is built from
$\Xi'$ then
\begin{equation}\label{eq:2}
\begin{split}
\frac 12 (c_1(\s _1) - c_1(\s _2)) = [\Sigma_3]^*
+ [\Sigma_4]^* - (p-1)[\Sigma_5]^* + \frac 12 ((x_1 - x_2)[\Sigma_6]^*  
& + (y_1 - y_2) [\Sigma_7]^* \\ & + (z_1 - z_2) [\Sigma_8]^*).
\end{split}
\end{equation}
The matrix formed by the top five rows of $\Phi^*$ is invertible over
$\Q$, therefore $[\Sigma_6]^*|_{W_{p,k,l,m}}$,
$[\Sigma_7]^*|_{W_{p,k,l,m}}$ and $[\Sigma_8]^*|_{W_{p,k,l,m}}$ are
linearly independent in $H^2 (X_{p,k,l,m}; \Z)$. This implies that
\[
\frac 12 ((x_1 - x_2) [\Sigma_6]^* + (y_1 - y_2)
[\Sigma_7]^* + (z_1 - z_2) [\Sigma_8]^*) 
\]
belongs to the image of $\varphi^*$ if and only if $x_1 - x_2=0$, $y_1
- y_2=0$, and $z_1 - z_2=0$. Thus if $\zeta_1$ and
$\zeta_2$ are both built from $\Xi$ or from $\Xi'$ then their surgery
presentations induce isomorphic spin$^c$ structures on $W_{p,k,l,m}$
if and only if $x_1 = x_2$, $y_1 = y_2$ and $z_1 = z_2$.

Let $c_i$ denote the $i$-th column of $\Phi^*$. The class  
\[
[\Sigma_3]^*+[\Sigma_4]^*-(p-1)[\Sigma_5]^*-[\Sigma_8]^*
\] 
can be expressed as $c_5-c_1-c_2$, and therefore its restriction
in $H^2(W_{p,k,l,m}; \Z)$ vanishes. Using this we see that 
Equation~\eqref{eq:2} implies
\begin{equation}\label{eq:3}
\frac 12 (c_1(\s _1) - c_1(\s _2))|_{W_{p,k,l,m}} = 
\frac 12 ((x_1 - x_2) [\Sigma_6]^* + (y_1 - y_2) [\Sigma_7]^* + 
(z_1 - z_2 +2) [\Sigma_8]^*)|_{W_{p,k,l,m}}.
\end{equation}
 Thus, by Equation~\eqref{eq:3} if $\zeta_1$ is obtained
by Legendrian surgery on $\Xi$ and $\zeta_2$ is obtained by Legendrian
surgery on $\Xi'$, then the surgery presentations of $\zeta_1$ and
$\zeta_2$ induce isomorphic spin$^c$ structures on $W_{p,k,l,m}$ if
and only if $x_1 = x_2$, $y_1 = y_2$ and $z_1 = z_2-2$.

The same argument given above works in the case of $M_{p,k,m}$. One
just needs to omit the knot $K_7$ from Figure~\ref{f:Mpklm} and work
with the analogously defined manifolds $W_{p,k,m}$ and
$X_{p,k,m}$. The new matrix $\Phi^*$ is obtained from the original
matrix $\Phi^*$ by simply dropping the seventh row. The remaining
computations are essentially the same, except one does not have terms
involving $y_1$, $y_2$ nor $[\Si_7]$. Similar considerations hold for
the cases of $M_{p,k,l}$ and $M_{p,k}$.
\end{proof}

\begin{prop}\label{p:donaldduck}
The number of isotopy classes of tight contact structures in 
$A(\Xi)$ is at least 
\[
(k-1)(l-1)m\quad\text{on}\quad
M_{p,k,l,m},
\]
\[
(k-1)m\quad\text{on}\quad
M_{p,k,m},
\]
\[
(k-1)(l-1)\quad\text{on}\quad
M_{p,k,l},
\]
and 
\[
(k-1)\quad\text{on}\quad
M_{p,k}.
\]
\end{prop}

\begin{proof}
In view of Lemma~\ref{l:goofy}, the number of different spin$^c$ 
structures induced on $W_{p,k,l,m}$ by the contact surgery diagrams of
Figure~\ref{f:Mpklm} gives a lower bound for the number of isotopy 
classes of
tight contact structures in $A(\Xi)$. Notice that $A(\Xi)$ can be decomposed
as $A\cup A'$, where $A$ contains the elements obtained by doing surgery on
$\Xi$, while $A'$ contains the ones obtained from $\Xi'$. By
Lemma~\ref{l:ziopaperone}(1) both $A$ and $A'$ contain $(k-1)(l-1)(m-1)$
elements distinguished by the induced spin$^c$ structures on
$W_{p,k,l,m}$. However, some elements may be contained both in $A$ and in
$A'$. In fact, by Lemma~\ref{l:ziopaperone}(2) for any contact structure
$\zeta$ in $A'$ there is a contact structure in $A$ inducing an isomorphic
spin$^c$ structure on $W_{p,k,l,m}$ unless ${\tt rot}(K_8)=-m$ in the surgery
diagram for $\zeta$. Since the number of contact surgery diagrams on
$M_{p,k,l,m}$ with ${\tt rot}(K_8)=-m$ giving tight contact structures
belonging to $A'$ is $(k-1)(l-1)$, there are at least
\[
(k-1)(l-1)(m-1) + (k-1)(l-1) = (k-1)(l-1)m
\]
nonisotopic tight contact structures in $A(\Xi)$. 

In the case of $M_{p,k,m}$, a similar argument gives the lower bound
$(k-1)m$. In the cases of $M_{p,k,l}$ and $M_{p,k}$, since there is no
knot $K_8$ every spin$^c$ structure induced by an element of $A'$ is
also induced by an element of $A$. Therefore, as a lower bound we just
get the number of elements of $A$, that is $(k-1)(l-1)$ in the case of
$M_{p,k,l}$ and $(k-1)$ in the case of $M_{p,k}$.
\end{proof}

\begin{cor}\label{c:AChi}
The number of isotopy classes of tight contact structures on the
3--manifold $M(-1; r_1, r_2, r_3)$ belonging to $A(\Xi)$ is at most
\[
(a_1^1-1)(a_1^2-1)a_1^3
\prod _{i=1}^3\prod_{j\geq 2}(a_j^i - 1)
\]
if $r_3\neq\frac1{a_0^3}$, and 
\[
\prod_{i=1}^2\prod_{j\geq 1}(a^i_j-1)
\]
if $r_3=\frac1{a_0^3}$. In the above formulae $a_j^i=2$ by convention
if $j>k_i$.
\end{cor}

\begin{proof}
The statement follows immediately from Proposition~\ref{p:donaldduck}
together with Theorem~\ref{t:fotight}.
\end{proof}

\begin{proof}[Proof of Theorem~\ref{t:main}]
The statement follows immediately combining Corollaries~\ref{c:ub1},
\ref{c:ub2}, \ref{c:Axi} and~\ref{c:AChi}.
\end{proof}

\end{document}